\documentclass[11pt]{amsart}
\usepackage[OT4]{fontenc}
\usepackage{graphicx, amsmath, amssymb, amscd, color}
\newtheorem{Theorem}{Theorem}
\newtheorem{Corollary}[Theorem]{Corollary}
\newtheorem{Lemma}[Theorem]{Lemma}
\newtheorem{Proposition}[Theorem]{Proposition}
\newtheorem*{example}{Example}
\newtheorem{prop}[Theorem]{Proposition}
\newtheorem{theorem}[Theorem]{Theorem}
\newtheorem{cor}[Theorem]{Corollary}
\newtheorem{lemma}[Theorem]{Lemma}
\newtheorem*{Assertion}{Assertion}
\theoremstyle{definition}
\newtheorem*{defn}{Definition}
\newtheorem*{definition}{Definition}
\newcommand{\norm}[1]{\left\Vert#1\right\Vert}

\newcommand{\set}[1]{\left\{#1\right\}}
\newcommand{\Real}{\mathbb R}
\newcommand{\R}{\mathbb R}
\newcommand{\Z}{\mathbb Z}

\newcommand{\sU}{\mathcal{U}}
\newcommand{\sV}{\mathcal{V}}
\newcommand{\sW}{\mathcal{W}}
\newcommand{\sH}{\mathcal{H}}
\newcommand{\di}{\partial_{\infty}}

\newcommand{\cal}{\mathcal}

\renewcommand{\epsilon}{\varepsilon}

 \DeclareMathOperator{\diam}{diam}
\DeclareMathOperator{\mesh}{mesh} \DeclareMathOperator{\Ind}{Ind}
 
\DeclareMathOperator{\dist}{dist} 
\DeclareMathOperator{\as}{asdim} \DeclareMathOperator{\asI}{asInd}
\DeclareMathOperator{\card}{Card} \DeclareMathOperator{\ind}{ind}

\begin{document}

\title{Asymptotic dimension}%
\author{G. Bell}%
\address{Mathematics \& Statistics\\ UNC Greensboro \\
383 Bryan Building \\
Greensboro NC, 27402 USA}%
\email{gcbell@uncg.edu}%
\author{A. Dranishnikov}
\address{Department of Mathematics\\ University of Florida \\
PO Box 118105\\ 358 Little Hall \\ Gainesville, FL 32611 USA}
\email{dranish@math.ufl.edu}
\thanks{The second author was partially
supported by NSF grant DMS-0604494}


\bibliographystyle{amsplain}

\begin{abstract}
The asymptotic dimension theory was founded by Gromov \cite{Gr93} in
the early 90s. In this paper we give a survey of its recent history
where we emphasize two of its features: an analogy with the
dimension theory of compact metric spaces and applications to the
theory of discrete groups.
\end{abstract}
\maketitle

\begin{quote} Counting dimensions we are definitely not counting
``things." \begin{flushright}{\it Yu. Manin},
\cite{Ma}\end{flushright}\end{quote}

%

\section{Introduction}

Dimension is a basic concept in mathematics. It came from Greek
geometry and made its way to all branches of modern mathematics and
science. Now dimension can be related to almost every mathematical
object.

In this survey we consider the asymptotic dimension of metric
spaces, in particular, discrete finitely generated groups taken with
word metrics. Asymptotic dimension theory bears a great deal of
resemblance to dimension theory in topology. In topology there are
three definitions of dimension: the covering dimension $dim$, the
large inductive dimension  $Ind$ and the small inductive dimension
$ind$. All three notions agree for separable metric spaces. In view
of the fact that $\Ind X=\dim X$ for all metrizable spaces whereas
$\ind X$ can be different, the dimension $\ind$ is less important
outside the class of separable metric spaces. For this reason, it
and its asymptotic analog will be avoided in this survey. We refer
to \cite{DrZ} for a definition of the asymptotic version of $\ind$.
The role of the asymptotic analog of $\Ind$ is limited to one
application in this stage of the theory. Thus, we will use the word
``dimension'' to mean the covering dimension and the term
``asymptotic dimension'' for the asymptotic analog of the covering
dimension.

Just as topological dimension is invariant under homeomorphisms,
its asymptotic analog is invariant under coarse isometries. Here a
{\it coarse isometry} is an isometry in the coarse category. {\em
The coarse category} can be described as follows. It is the
category whose objects are metric spaces $(X,d_X)$ and whose
morphisms are (not necessarily continuous) maps $f:(X,d_X)\to
(Y,d_Y)$ that are {\it metrically proper}, (i.e., the preimage of
every bounded set is bounded), and are {\it uniformly expansive},
(i.e., there is a function $\rho:\mathbb{R}\to \mathbb{R}$ such
that $d_Y(f(x),f(x'))\le \rho(d_X(x,x'))$). Two morphisms are said
to be equivalent if they are in finite distance from each other,
i.e., there is a constant $D>0$ such that $d_Y(f(x),g(x))<D$ for
all $x\in X$. The typical example of a coarse isometry is the
inclusion of the integers into the reals $\mathbb{Z}\to
\mathbb{R}$ supplied with the standard metric $d(x,y)=|x-y|$.
Also, we note that every bounded metric space is coarsely
equivalent to a point. For more details, see \S3.2 and \S12.

Asymptotic dimension was introduced by Gromov as an invariant of a
finitely generated group \cite{Gr93}. A (finite) symmetric
generating set $S$ on a group $\Gamma$ defines the word metric $d_S$
by taking the distance $d_S(x,y)$ to be the length of the shortest
presentation of $x^{-1}y$ in the alphabet $S$. Such a metric is
left-invariant. For a finitely generated group $\Gamma$, any choice
of generating sets gives rise to coarsely equivalent metric spaces
and hence, the asymptotic dimension $\as\Gamma$ is a group
invariant. Significant attention to asymptotic dimension was brought
by a theorem of Guoliang Yu \cite{Yu1} that proves the Novikov
higher signature conjecture for manifolds whose fundamental group
has finite asymptotic dimension. Later the integral Novikov
conjecture and some of its relatives were proved for groups with
finite asymptotic dimension
\cite{Ba,CG,Dr1,Dr3,DFW}. 
Unfortunately not all finitely presented groups have finite
asymptotic dimension. For example, Thompson's group $F$ has infinite
asymptotic dimension since it contains ${\mathbb Z}^n$ for all $n$.
Many relatives of the Novikov conjecture (like the Borel conjecture)
are of great interest for groups with finite cohomological
dimension. This makes the problem similar to the old Alexandroff
Problem from dimension theory: \emph{Does the cohomological
dimension of a space coincide with its covering dimension?} In
particular, this makes the asymptotic dimension theory a more
attractive subject.

This paper consists of two parts. In the first part we survey the
development of asymptotic dimension theory of metric spaces. In the
second part we discuss applications to groups and we present some
computations of the asymptotic dimension and some finite
dimensionality results. Our survey is in no way complete. We don't
discuss the many now-existing generalizations of asymptotic
dimension, nor do we present any applications of asymptotic
dimension to coarse geometry and to the Novikov-type conjectures; we
omit entirely cohomological asymptotic dimension theory.

The authors would like to thank Yehuda Shalom for some valuable
remarks.

\

\

\

{\bf I. ASYMPTOTIC DIMENSION OF SPACES}

\

\section{\bf Basic facts of classical dimension theory}

\subsection{Definitions} The definition of the covering dimension $\dim$ of a
topological space $X$ is due to Lebesgue: {\it $\dim X\le n$ if and
only if for every open cover $\sV$ of $X$ there is a cover $\sU$ of
$X$ refining $\sV$ with multiplicity $\le n+1.$}

We will use the words \textit{order} and \textit{multiplicity} of a
cover interchangeably to mean the largest number of elements of the
cover meeting any point of the the space. Given a cover $\sV$ of a
topological space, we say that the cover $\sU$ \textit{refines}
$\sV$ if every $U\in\sU$ is contained in some element $V\in\sV.$

As usual, we define $\dim X=n$ if it is true that $\dim X\le n$ but
it is not true that $\dim X\le n-1.$

Recall that a closed subset $C$ of a topological space $X$ is a {\it
separator\/} between disjoint subsets $A, B\subset X$ if $X\setminus
C=U\cup V$, where $U$ and $V$ are open subsets in $X$, $U\cap
V=\emptyset$, $A\subset U$, $B\subset V$. A closed subset $C$ of a
topological space $X$ is a {\it cut} between disjoint subsets $A$
and $B\subset X$ if every continuum (compact connected space)
$T\subset X$ that intersects both $A$ and $B$ also intersects $C$.
The definition of large inductive dimension is due to Brouwer and
Poincare: {\it $Ind X\le n$ if for every pair of closed disjoint
sets $A$ and $B\subset X$ there is a separator $C$ with $Ind C\le
n-1$ where $Ind X=-1$ if and only if $X=\emptyset$.} By replacing
separators with cuts we obtain the definition of Brouwer's {\it
Dimensiongrad}, $Dg(X)$. For compact Hausdorff spaces there are the
inequalities:
$$ \dim X\le Dg(X)\le Ind X.$$

The definition of the small inductive dimension is due to Menger and
Urysohn: $ind X=-1$ if and only if $X=\emptyset$; $ind X\le n$ if
for every point $x\in X$ and every neighborhood $U$ of $x$ there is
a smaller neighborhood $V$, $x\in V\subset U$ with $ind(\partial
V)\le n-1$.

\begin{Theorem}
For a compact metric space $X$ the following conditions are
equivalent.

(1) $dim X\le n$;

(2) $Ind X\le n$;

(3) $ind X\le n$;

(4) $Dg(X)\le n$;

(5) Every continuous map $f:A\to S^n$ of a closed subset $A\subset
X$ to the $n$-sphere has a continuous extension $\bar f:X\to S^n$;

(6) For every $\epsilon>0$ there is an $\epsilon$-map $\phi:X\to
K^n$ to an $n$-dimensional polyhedron;

(7) $X$ is the limit of an inverse sequence,
$X=\lim_{\leftarrow}K_i$, of $n$-dimensional polyhedra $K_i$;

(8) For every $\epsilon>0$ there is an $\epsilon$-cover ${\sU}$ of
$X$ which can be decomposed into $n+1$ disjoint families
${\sU}={\sU}^0\cup\dots\cup{\sU}^n$;

(9) $X$ can be presented as the union of $n+1$ $0$-dimensional
subsets;

(10) $X$ admits a light map $f:X\to I^n$ onto the $n$ cube.
\end{Theorem}

The equivalence of (4) to all other was proven in \cite{FLS}. All
other equivalences are well-known facts and can be found in any
textbook (see, for example \cite{En95}). One can extend these
equivalences (probably with the exception of (4)) to all separable
metric spaces by replacing (where it is appropriate) an arbitrary
$\epsilon$ by an arbitrary open cover of $X$.

We now recall all the terminology used in the above theorem: A map
$\phi : X\to Y$ of a metric space is an {\it $\epsilon$-map} if
$\diam(\phi^{-1}(y))<\epsilon$ for all $y\in Y$; A cover $\cal U$ is
an {\em $\epsilon$-cover} if $\diam U<\epsilon$ for all $U\in\cal
U$; A family $\sV$ of subsets of a space $X$ is {\it disjoint} if
$V\cap V'=\emptyset$ for all $V,V'\in\sV$, $V\ne V'$; A map $f:X\to
Y$ is called {\it light} if $\dim f^{-1}(y)=0$ for all $y\in Y$.

Note that for every closed subset $Y\subset X$, $\dim Y\le\dim X$.
This is also true for every subset of a metric space and is not
generally true for every subset of compact (non-metrizable) spaces
\cite{En95}.
\subsection{Union Theorems}
\begin{Theorem} Let $X$ be a separable metric space, then:

1. $\dim X=\max\{\dim A,\dim(X\setminus A)\}$ for every closed
subset $A\subset X$.

2.(Countable Union Theorem) Let $X=\cup_{i=1}^{\infty} X_i$ be the
countable union of closed subsets, then $\dim X=\sup_i\{\dim X_i\}$.

3. (Menger-Urysohn sum Formula) $\dim(X\cup Y)\le\dim X+\dim Y+1$
for arbitrary sets.
\end{Theorem}

\subsection{Dimension and mappings}
\begin{Theorem}[Hurewicz Mapping Theorem 1]\label{dim Hurewicz
Theorem} Let $f:X\to Y$ be a map between compact spaces. Then
$$
\dim X\le\dim Y+\dim f
$$
where $\dim f=\sup\{\dim f^{-1}(y)\mid y\in Y\}$.
\end{Theorem}

The Hurewicz Mapping Theorem 1 implies that a light map of
compacta cannot lower the dimension, in particular a map with
finite preimages cannot lower the dimension:
\begin{Theorem}[Hurewicz Mapping Theorem 2]
Let $f:X\to Y$ be a map between compact spaces with $|f^{-1}(y)|\le
k$. Then
$$
\dim X\le\dim Y+k-1
$$
\end{Theorem}

We note that these theorems hold for all metric spaces and closed
mappings.

\subsection{Dimension of the product} Most of the theorems on the
dimension of product are proven by cohomological approach.
\begin{Theorem} {\it For all normal spaces $\dim(X\times Y)\le\dim X+\dim
Y$.}
\end{Theorem}

Note that when one of the factors is compact, this theorem follows
from the Hurewicz Mapping Theorem 1.

\begin{Theorem}[Morita Theorem] For all normal spaces $X$,
$\dim(X\times[0,1])=\dim X+1$ .
\end{Theorem}
The equality does not hold for 2-dimensional subsets
$X\subset{\Real}^3$ if one replaces the interval $[0,1]$ by a
continuum \cite{DRS2}, although the following does hold:
\begin{Theorem} {\it For a compact space $X$ and a continuum $Y$ there
is the inequality: $\dim(X\times Y)\ge\dim X+1$.}
\end{Theorem}

\begin{Theorem} {\it For a compact metric space $X$ either $\dim (X\times
X)=2\dim X$ or $\dim(X\times X)=2\dim X-1$.}
\end{Theorem}
This theorem divides all compacta in two classes: Type I and Type
II. There is a theory of Bockstein that allows one to determine the
type of a compactum by means of its cohomological dimension with
respect to different coefficient groups \cite{Dr-Cohom-Compact}.
This theory allows one to compute the dimension of the product. The
Bockstein theory combined with the Realization Theorem of
\cite{DR-Mapping-intersection} implies
\begin{Theorem} {\it For any natural numbers $m$,$n$, and any $k$ with
$\max\{m,n\}<k\le m+n$, there are compact metric spaces $X$ and $Y$
with the dimensions $\dim X=m$, $\dim Y=n$, $\dim(X\times Y)=k$.}
\end{Theorem}
Cohomological dimension theory allows one to improve the union
theorem and the Hurewicz Mapping Theorem 1:
\begin{Theorem} {\it Suppose that $X\cup Y$ is a compactum, then

(1) $\dim(X\cup Y)\le\dim(X\times Y)+1$ if $X\cup Y$ is of the first
type;

(2) $\dim(X\cup Y)\le\dim(X\times Y)+2$ if $X\cup Y$ is of the
second type.}
\end{Theorem}
There are examples where the inequality (2) is sharp.
\begin{Theorem} {\it Let $f:X\to Y$ be a map between compacta. Then

(1) $\dim X\le\sup\{\dim(Y\times f^{-1}(y))\mid y\in Y\}$ if $X$ is
of the first type;

(2) $\dim X\le\sup\{\dim(Y\times f^{-1}(y))\mid y\in Y\}+1$ if $X$
is of the second type.}
\end{Theorem}
We note that in both theorems, statement (2) is not always an
improvement on the classic results.

\subsection{Embedding theorems}
We denote the Stone-\v Cech compactification of a space $X$ by
$\beta X$.
\begin{Theorem} $\dim X=\dim\beta X$.
\end{Theorem}
This together with the Marde\v si\' c Factorization theorem implies
the existence of compact metric spaces of dimension $n$ that contain
a topological copy of every separable metric space of dimension $n$.
Such compacta are called {\it universal}.
\begin{Theorem}[Marde\v si\' c Factorization Theorem]
{\it For every continuous map $f:X\to Y$ of a compact Hausdorff
space to a metric space there are maps $g:X\to X'$ and $f':X'\to Y$
such that $X'$ is metrizable and $f=f'\circ g$.}
\end{Theorem}
There are ``nice'' universal compacta, namely Menger compacta
$\mu^n$, which are characterized by the following \cite{Bestvina}:
\begin{Theorem}[Bestvina Criterion]{\it A compact metric space $X$ is
homeomorphic to the Menger compactum $\mu^n$ if and only if it is
$n$-dimensional, $n-1$-connected and locally $n-1$-connected, and it
satisfies the disjoint $n$-disc property, $DDP^n$, i.e., every two
maps of the $n$-disc $f,g:D^n\to X$ can be approximated by maps with
disjoint images.}
\end{Theorem}
In non-compact case there are universal spaces $\nu^n$ in dimension
$n$ called N\"obeling spaces, defined as the set of all points in
${\Real}^{2n+1}$ with at most $n$ rational coordinates. An analogous
topological characterization of the N\"obeling space $\nu^n$ was
given by Nagorko \cite{Nag} (see also \cite{Levin} for a different
treatment). Since by its construction $\nu^n$ is a subset of
${\Real}^{2n+1}$, we obtain
\begin{Theorem} [N\"obeling-Pontryagin Theorem]\label{NP}
{\it Every $n$-dimensional compact metric space can be embedded into
${\mathbb R}^{2n+1}$.}
\end{Theorem}
\begin{Theorem}\cite{DRS1,Sp}
{\it Every compact $n$-dimensional metric space of the second type
admits an embedding into ${\mathbb R}^{2n}$}.
\end{Theorem}
In these two theorems as well as in the case of a Menger space,
every map of an $n$-dimensional compactum to $\mu^n$, ${\mathbb
R}^{2n+1}$ (Theorem 15) or to ${\mathbb R}^{2n}$ (Theorem 16) can be
approximated by embeddings.

We recall that a {\it dendrite} is an 1-dimensional, 1-connected,
locally 1-connected compact metric space.
\begin{Theorem}\cite{Bowers} \label{Thm17} {\it Every compact metric space $X$ of dimension $\dim
X\le n$ can be embedded into the product $\prod_{i=1}^n
T_i\times[0,1]$ of $n$ dendrites and an interval.}
\end{Theorem}
A similar result was obtained in \cite{St}.
\subsection{Infinite dimensional spaces}
A space $X$ is called \emph{strongly infinite dimensional} if it
admits an essential map onto the Hilbert cube
$I^{\infty}=\prod_{i=1}^{\infty}[0,1]$. We recall that a map $f:X\to
I^{\infty}$ is called \emph{essential} if every projection onto a
finite-dimensional face $\pi:I^{\infty}\to I^n$ is essential. A map
$g:X\to I^n$ is essential if it cannot be deformed to a map
$g':X\to\partial I^n$ by a deformation fixed on $g^{-1}(\partial
I^n)$. An infinite dimensional space is called \emph{weakly infinite
dimensional} if it is not strongly infinite dimensional. A space $X$
is called \emph{countable dimensional} if it can be presented as a
countable union of 0-dimensional subspaces. A space $X$ has
\emph{Property C} if for every sequence of open covers
${\sV}_1,\dots,{\sV}_k,\dots$ of $X$ there are disjoint families of
open subsets ${\sU}_1,\dots,{\sU}_k,\dots$ such that each $\sU_i$ is
a refinement of ${\sV}_i$ and $\cup_i{\sU}_i$ is a cover of $X$.
\begin{Theorem}
Among compact metric spaces there are inclusions:
$$
\{\text{Countable dim}\}\subset\{\text{Property C}\}\subset\{\text
{Weakly infinite dim}\}.$$
\end{Theorem}
Due to examples of R. Pol \cite{Pol} and P. Borst \cite{Borst}, both
inclusions are strict. Most of the theorems for finite dimensional
compacta can be extended to compacta with Property C, but usually
not to the class of strongly infinite dimensional compacta. For
example, the Alexandroff Theorem stating that the covering dimension
$\dim X$ coincides with the cohomological dimension $dim_{\mathbb
Z}X=\sup\{n\mid \check H^n(X,A)\ne 0, A\subset_{Cl}X\}$ for finite
dimensional compacta can be extended to compacta with Property C
\cite{An}, but cannot be extended to strongly infinite dimensional
compacta \cite{Dr-Cohom-Compact}.

\section{Definitions of asdim}

\subsection{Equivalent definitions}
\begin{defn} Let $X$ be a metric space.  We say that the \emph{asymptotic
dimension} of $X$ does not exceed $n$ and write $\as X\le n$
provided for every uniformly bounded open cover $\sV$ of $X$ there
is a uniformly bounded open cover $\sU$ of $X$ of multiplicity $\le
n+1$ so that $\sV$ refines $\sU.$ We write $\as X=n$ if it is true
that $\as X\le n$ and $\as X\not\le n-1$.
\end{defn}
Note that the asymptotic dimension can be viewed as somehow dual to
Lebesgue covering dimension.

Often we will need to consider very large positive constants and we
remind ourselves that they are large by writing $r<\infty$ instead
of $r>0.$  On the other hand, writing $\epsilon>0$ is supposed to
mean that $\epsilon$ is a small positive constant.

Let $r<\infty$ be given and let $X$ be a metric space. We will say
that a family $\sU$ of subsets of $X$ is $r$-{\em disjoint} if
$d(U,U')>r$ for every $U\neq U'$ in $\sU.$ Here, $d(U,U')$ is
defined to be $\inf\{d(x,x')\mid x\in U, x'\in U'\}.$  The $r$-{\em
multiplicity} of a family $\sU$ of subsets of $X$ is defined to be
the largest $n$ so that no ball $B_r(x)\subset X$ meets more than
$n$ of the sets from $\sU$; more succinctly, the $r$-multiplicity of
$\sU$ is $\sup_{x\in X}\card\{U\in\sU\mid U\cap
B_r(x)\neq\emptyset\}$. Recall that the {\em Lebesgue number} of a
cover $\sU$ of $X$ is the largest number $\lambda$ so that if
$A\subset X$ and $\diam(A)\le\lambda$ then there is some $U\in\sU$
so that $A\subset U.$

Let $K$ be a countable simplicial complex.  There are two natural
metrics one can place on is geometric realization $|K|$: the uniform
metric and the geodesic metric. We wish to consider the uniform
metric on $|K|.$ This is defined by embedding $K$ into $\ell^2$ by
mapping each vertex to an element of an orthonormal basis for
$\ell^2$ and giving it the metric it inherits as a subspace.  A map
$\varphi:X\to Y$ between metric spaces is uniformly cobounded if for
every $R>0$, $\diam(\varphi^{-1}(B_R(y)))$ is uniformly bounded.

\begin{Theorem} \label{equiv asdim} Let $X$ be a metric space.  The following conditions
are equivalent.
\begin{enumerate}
  \item $\as X\le n$;
  \item for every $r<\infty$ there exist $r$-disjoint families
  $\sU^0,\ldots, \sU^n$ of uniformly bounded subsets of $X$
  such that $\cup_i\sU^i$ is a cover of $X$;
  \item for every $d<\infty$ there exists a uniformly bounded
  cover $\sV$ of $X$ with $d$-multiplicity $\le n+1$;
  \item for every $\lambda<\infty$ there is a uniformly bounded
  cover $\sW$ of $X$ with Lebesgue number $>\lambda$ and multiplicity
  $\le n+1$; and
  \item for every $\epsilon>0$ there is a uniformly cobounded,
  $\epsilon$-Lipschitz map $\varphi:X\to K$ to a uniform simplicial
  complex of dimension $n.$
\end{enumerate}
\end{Theorem}
The proof can be found in \cite{BD4}.

We conclude this section with an example.

\begin{example} \label{asdimT=1} $\as T\le 1$ for all trees $T$ in the edge-length metric.
\end{example}

\begin{proof} Fix some vertex $x_0$ to be the root of the tree.  Let
$r<\infty$ be given and take concentric annuli centered at $x_0$ of
thickness $r$ as follows: $A_k=\{x\in T\mid
d(x,x_0)\in[kr,(k+1)r)\}.$  Although alternating the annuli (odd k,
even k) yields $r$-disjoint sets, these sets clearly do not have
uniformly bounded diameter.  We have to further subdivide each
annulus.

Fix $k>1.$ Define $x\sim y$ in $A_k$ if the geodesics $[x_0,x]$ and
$[x_0,y]$ in $T$ contain the same point $z$ with
$d(x_0,z)=r(k-\frac12).$ Clearly in a tree this forms an equivalence
relation.  The equivalence classes are $3r$ bounded and elements
from distinct classes are at least $r$ apart.  So, define $\sU$ to
be equivalence classes corresponding to even $k$ (along with $A_0$
itself) and $\sV$ to be equivalence classes corresponding to odd
$k.$  These two families cover $T$ and consist of uniformly bounded,
$r$-disjoint sets.  Thus, $\as T\le 1.$
\end{proof}

\subsection{Large-scale invariance of asdim}

Let $X$ and $Y$ be metric spaces. A map $f:X\to Y$ between metric
spaces is a {\em coarse embedding} if there exist non-decreasing
functions $\rho_1$ and $\rho_2$,
$\rho_i:\overline{\mathbb{R}}_+\to\overline{\mathbb{R}}_+$ such that
$\rho_i\to\infty$ and for every $x,x'\in X$
\[\rho_1(d_X(x,x'))\le d_Y(f(x),f(x'))\le
\rho_2(d_X(x,x'))\text{.}\]  Such a map is often called a {\em
coarsely uniform embedding} or just a {\em uniform embedding}.  The
metric spaces $X$ and $Y$ are {\em coarsely equivalent} if there is
a coarse embedding $f:X\to Y$ so that there is some $R$ such that
$Y\subset N_R(f(X)).$ Equivalently, to metric spaces are coarsely
equivalent if there exist coarse embeddings $f:X\to Y$ and $g:Y\to
X$, whose compositions (in both ways) are $K$-close to the identity,
for some $K>0$.

A coarse embedding $f:X\to Y$ is a {\em quasi-isometric embedding}
if it admits linear $\rho_i.$ The two spaces $X$ and $Y$ are
quasi-isometric if there is a quasi-isometry $f:X\to Y$ and a
constant $C$ so that $Y \subset N_C(f(X)).$

A metric space $(X,d)$ is called \emph{geodesic} if for every two
points $x,y\in X$ there is an isometric embedding of the interval
$\xi:[0,a]\to X$ with $a=d(x,y)$, $\xi(0)=x$, and $\xi(a)=y$. We
note that a coarse embedding of a geodesic metric space always
admit a linear function $\rho_2$. This implies that a coarse
equivalence between geodesic metric spaces is in fact a
quasi-isometry.

A metric space $(X,d)$ is called \emph{$c$-discrete} if $d(x,x')\ge
c$ for all $x,x'\in X$, $x\ne x'$.

\begin{Proposition}\label{1-discrete}
Every metric space $(X,d)$ is coarsely equivalent to a 1-discrete
metric space.
\end{Proposition}
\begin{proof}
Let $S\subset X$ be a maximal 1-discrete subset. Then the inclusion
is a coarse equivalence.
\end{proof}

\begin{Corollary}
Every geodesic metric space is quasi-isometric to a connected graph.
\end{Corollary}
\begin{proof}
Let $S\subset X$ be a maximal 1-discrete subset in $X$. Connect all
points $s\neq s'$ in $S$ with $d(s,s')\le 4$ by an edge. Let $\rho$
be a simplicial metric on the graph (so that every edge has length
one). Clearly, $\rho(s,s')\ge \frac{1}{3}d_X(s's')$. We show that
$\rho(s,s')\le d_X(s's')+2$. Let $x_0,x_1,\dots, x_n$ be points on a
geodesic segment joining $s$ and $s'$ such that $s=x_0$, $s'=x_n$,
and $d_X(x_{i-1},x_{i})=1$ for $i<n$. Let $s_i\in S$ be such that
$d_X(s_i,x_i)\le 1$. Then $s_i$ is joined by an edge with $s_{i+1}$
for all $i$. Hence $\rho(s,s')\le n+2\le d_X(s,s')+2$.
\end{proof}

\begin{Proposition} Let $f:X\to Y$ be a coarse equivalence.  Then
$\as X=\as Y.$
\end{Proposition}

\begin{proof} If $\sU^0,\ldots,\sU^n$ are $r$-disjoint, $D$-bounded
families covering $X$ then the families $f(\sU^i)$ are
$\rho_1(r)$-disjoint and $\rho_2(D)$-bounded.  Since $N_R(f(X))$
contains $Y$ we see that taking families $N_R(f(\sU^i))$ will cover
$Y$ and be $(2R+\rho_2(D))$-bounded and $(\rho_1(r)-2R)$-disjoint.
Since $\rho_i\to\infty,$ $r$ can be chosen large enough for
$\rho_1(r)-2R$ to be as large as one likes.  Therefore, $\as Y\le
\as X.$

The same proof applied to a coarse inverse for $f$ proves that $\as
X\le \as Y.$
\end{proof}

\begin{example}  $\as\mathbb{R}=\as\mathbb{Z}=1.$
\end{example}

\begin{Proposition} Let $X$ be a metric space and $Y\subset X.$
Then $\as Y\le \as X.$
\end{Proposition}

\begin{proof} Let $R<\infty$ be given and take a cover $\sU$ of $X$
by uniformly bounded sets with $R$-multiplicity $\le n+1.$  Clearly
the restriction of this cover to $Y$ yields a cover whose elements
have uniformly bounded diameter and at most $n+1$ of them can meet
any ball of radius $R$ in $Y.$  Thus, $\as Y\le \as X.$
\end{proof}

\section{Union theorems}
Before proceeding further, we need to establish a basic result: a
union theorem for asymptotic dimension. It should be noted that here
asymptotic dimension varies slightly from covering dimension. For
example, the finite union theorem for covering dimension says
$\dim(X\cup Y)\le \dim X+\dim Y+1,$ and that inequality is sharp.
Also, the standard countable union theorem for covering dimension is
$\dim(\cup_i C_i)\le \max_{i}\{\dim C_i\}$ where the $C_i$ are
closed subsets of $X.$ Notice that there can be no direct analog of
this theorem for asymptotic dimension since every finitely generated
group is a countable set of points, and as we shall see, finitely
generated groups can have arbitrary (even infinite) asymptotic
dimension.

Let $\sU$ and $\sV$ be families of subsets of $X$.  Define the
$r$-saturated union of $\sV$ with $\sU$ by
\[\sV\cup_r\sU=\{N_r(V;\sU)\mid V\in\sV\}\cup \{U\in\sU\mid d(U,\sV)>r\},\]
where $N_r(V;\sU)=V\cup\bigcup_{d(U,V)\le r} U.$

It is easy to verify the following proposition.

\begin{Proposition} Let $\sU$ be an $r$-disjoint, $R$-bounded family
of subsets of $X$ with $R\ge r.$  Let $\sV$ be a $5R$-disjoint, $D$
bounded family of subsets of $X.$  Then $\sV\cup_r\sU$ is
$r$-disjoint and $(D+2(r+R))$-bounded.\qed
\end{Proposition}


Let $X$ be a metric space.  We will say that the family
$\{X_\alpha\}$ of subsets of $X$ satisfies the inequality $\as
X_\alpha\le n$ {\em uniformly} if for every $r<\infty$ one can find
a constant $R$ so that for every $\alpha$ there exist $r$-disjoint
families $\sU^0_\alpha,\ldots,\sU^n_\alpha$ of $R$-bounded subsets
of $X_\alpha$ covering $X_\alpha.$ A typical example of such a
family is a family of isometric subsets of a metric space. Another
example is any family containing finitely many sets.

\begin{Theorem}[Union Theorem] Let $X=\cup_\alpha X_\alpha$ be a metric space where
the family $\{X_\alpha\}$ satisfies the inequality $\as X_\alpha\le
n$ uniformly.  Suppose further that for every $r$ there is a
$Y_r\subset X$ with $\as Y_r\le n$ so that
$d(X_\alpha-Y_r,X_{\alpha'}-Y_r)\ge r$ whenever $X_\alpha\neq
X_{\alpha'}.$  Then $\as X\le n.$
\end{Theorem}

Before proving this theorem, we state a corollary: the finite union
theorem for asymptotic dimension.

\begin{Corollary}[Finite Union Theorem] Let $X$ be a metric space with $A,B\subset X$.
Then $\as (A\cup B)\le\max\{\as A, \as B\}.$
\end{Corollary}

\begin{proof}[Proof of Corollary.] Apply the union theorem to the
family $\{A,B\}$ with $B=Y_r$ for every $r.$
\end{proof}

\begin{proof}[Proof of Union Theorem.] Let $r<\infty$ be given and
take $r$-disjoint, $R$-bounded families $\sU^i_\alpha$ $(i=0,\ldots,
n)$ of subsets of $X_\alpha$ so that $\cup_i\sU^i_\alpha$ covers
$X_\alpha.$  We may assume $R\ge r.$ Take $Y=Y_{5R}$ as in the
statement of the theorem and cover $Y$ by families
$\sV^0,\ldots,\sV^n$ that are $D$-bounded and $5R$-disjoint. Let
$\bar\sU^i_\alpha$ denote the restriction of $\sU^i_\alpha$ to the
set $X_\alpha-Y.$  For each $i,$ take
$\sW^i_\alpha=\sV^i\cup_r\bar\sU^i_\alpha.$ By the proposition
$\sW^i$ consists of uniformly bounded sets and is $r$-disjoint.
Finally, put $\sW^i=\{W\in \sW^i_\alpha\mid \alpha\}$. Observe that
each $\sW^i$ is $r$-disjoint and uniformly bounded. Also, it is easy
to check that $\cup_i\sW^i$ covers $X.$
\end{proof}

\section{Connection with covering dimension - Higson corona}
Let $\varphi\colon X\to\mathbb{R}$ be a function defined on a metric
space $X$. For every $x\in X$ and every $r>0$ let
$V_r(x)=\sup\{|\varphi(y)-\varphi(x)|\colon y\in N_r(x)\}$. A
function $\varphi$ is called {\it slowly oscillating\/} if for every
$r>0$ we have $V_r(x)\to0$ as $x\to\infty$. (This means that for
every $\varepsilon>0$ there exists a compact subspace $K\subset X$
such that $|V_r(x)|<\varepsilon$ for all $x\in X\setminus K$). Let
$\bar X$ be the compactification of $X$ that corresponds to the
family of all continuous bounded slowly oscillating functions. The
{\it Higson corona\/} of $X$ is the remainder $\nu X=\bar X\setminus
X$ of this compactification.

It is known that the Higson corona is a functor from the category of
proper metric space and coarse maps into the category of compact
Hausdorff spaces. In particular, if $X\subset Y$, then $\nu X\subset
\nu Y$.

For any subset $A$ of $X$ we denote by $A'$ its trace on $\nu X$,
i.e. the intersection of the closure of $A$ in $\bar X$ with $\nu
X$. Obviously, the set $A'$ coincides with the Higson corona $\nu
A$.

Dranishnikov, Keesling and Uspenskij \cite{DKU} proved the
inequality
$$
\dim\nu X\le\as X,
$$
for any proper metric space $X$. It was shown there that $\dim\nu
X\ge\as X$ for a large class of spaces, in particular for
$X=\mathbb{R}^n$. Later Dranishnikov proved \cite{Dr1} that the
equality $\dim\nu X=\as X$ holds provided $\as X<\infty$. The
question of whether there is a metric space $X$ with $\as X=\infty$
and $\dim\nu X<\infty$ is still open.

\section{Inductive approach}
The notion of asymptotic inductive dimension $\Ind$ was introduced
in \cite{Dr2}.

Let $X$ be a proper metric space. A subset $W\subset X$ is called an
{\it asymptotic neighborhood\/} of a subset $A\subset X$ if
$\lim_{r\to\infty}d(A\setminus N_r(x_0),X\setminus W)=\infty$. Two
sets $A,B$ in a metric space are {\it asymptotically disjoint\/} if
$\lim _{r\to\infty}d(A\setminus N_r(x_0),B\setminus
N_r(x_0))=\infty$. In other words, two sets are asymptotically
disjoint if the traces $A'$, $B'$ on $\nu X$  are disjoint.

A subset $C$ of a metric space $X$ is an  {\it asymptotic
separator\/}  between asymptotically disjoint subsets
$A_1,A_2\subset X$ if the trace $C'$ is a separator in $\nu X$
between $A_1'$ and $A_2'$.

We recall the definition of the asymptotic Dimensiongrad in the
sense of Brouwer, $asDg$, from \cite{DrZ}.

Let $X$ be a metric space and $\lambda>0$. A finite sequence
$x_1,\dots,x_k$ in $X$ is a {\it $\lambda$-chain\/} between subsets
$A_1,A_2\subset X$ if $x_1\in A_1$, $x_k\in A_2$ and
$d(x_i,x_{i+1})<\lambda$ for every $i=1,\dots, k-1$. We say that a
subset $C$ of a metric space $X$ is an  {\it asymptotic cut\/}
between the asymptotically disjoint subsets $A_1,A_2\subset X$ if
for every $D>0$ there is a $\lambda>0$ such that every
$\lambda$-chain between $A_1$ and $A_2$ intersects $N_D(C)$.

By definition, $\asI X=asDg(X)=-1$ if and only if $X$ is bounded.
Suppose we have defined the class of all proper metric spaces $Y$
with $\asI Y\le n-1$ (respectively with $asDg(Y)\le n-1$). Then
$\asI X\le n$ (respectively  $asDg(Y)\le n$) if and only if for
every asymptotically disjoint subsets $A_1,A_2\subset X$ there
exists an asymptotic separator (respectively asymptotic cut) $C$
between $A_1$ and $A_2$ with $\asI C\le n-1$ (respectively $asDg
(C)\le n-1$). The dimension functions $\asI$ and $asDg$ are called
the {\em asymptotic inductive dimension\/} and {\em asymptotic
Brouwer inductive dimension\/} respectively.

\begin{Proposition} $asDg(X)\le\asI X$, for
every $X$. \end{Proposition} It is unknown if $\asI=asDg$ for proper
metric spaces.

\begin{Theorem}\label{asInd le asdim}
For all proper metric spaces $X$ with $0<\as X<\infty$ we have
$$
\as X=\asI X.
$$
\end{Theorem}

This theorem is a very important step in the  proof of the exact
formula of the asymptotic dimension of the free product $\as A\ast
B$ of groups \cite{BDK}, see section 3.4.

Notice that there is a small problem with coincidence of $\as$ and
$\Ind$ in dimension $0.$  This leads to philosophical discussions of
whether bounded metric spaces should be defined to have $\as=-1$ or
$0.$ Observe that in the world of finitely generated groups, $\as
\Gamma=0$ if and only if $\Gamma$ is finite. On the other hand,
there are metric spaces, for instance $2^n\subset\mathbb{R}$ that
are unbounded yet have asymptotic dimension $0.$

\subsection{Proof that asInd X$\le$ asdim X}

Let $N$ be a simplicial complex. For a subset $A\subset N$ we denote
by $st(A)$ the {\em star neighborhood of $A$}, i.e., the union of
all simplices in $N$ that have nonempty intersection with $A$. Note
that $st(A)$ is a subcomplex of $N$. The {\em canonical regular
neighborhood $W$} of a subcomplex $K\subset N$ is the star
neighborhood of $K$ in the second regular barycentric subdivision of
$N$. The {\em regular neighborhood} $W$ is the mapping cylinder
neighborhood $M_\phi$ of a simplicial map $\phi:\partial W\to K$,
with respect to the second barycentric subdivision of $N$.

For a function $f:X\to\R$ we denote its ``coarse derivative" at a
point $x$ by \[\nabla f(x)=\diam(f(B_1(x)))\hbox{.}\]

\begin{prop} \label{Prop 1}Let $f:N\to [-1,1]$ be a map defined on an $n$-dimensional
uniform simplicial complex that is extendible to the Higson corona.
Suppose that $f^{-1}(0)\setminus A$ is an $(n-1)$-dimensional
subcomplex of $N$ for some set $A\subset f^{-1}(0)$. Then there is a
function $\tilde f:N\to[-1,1]$ with $|f(x)-\tilde f(x)|\le 2\nabla
f(x)$ such that $\tilde f^{-1}(0)$ is an $(n-1)$-dimensional
subcomplex of $N$ and $\tilde f^{-1}(0)\subset f^{-1}(0)\cup W$,
where $W$ is the regular neighborhood of $st(A)$.
\end{prop}

\begin{proof} Let $q:\partial W\times [0,1]\coprod st(A)\to W$ be the
quotient map from the definition of the mapping cylinder. Define
\[\tilde f_1(x)=\left\{
   \begin{array}{ll}
     0, & \hbox{if $x\in st(A)$;} \\
     f(x), & \hbox{if $n\in N\setminus W$;} \\
     (1-t)f(y), & \hbox{if $x=q(y,t)$.}
   \end{array}
 \right.
\]
For every $x\in W$ there is a $z_x\in A$ with $d(x,z_x)<2$. Note
that $|\tilde f_1(x)-f(x)|=|f(x)|=|f(x)-f(z_x)|\le \nabla f(x)$ when
$x\in st(A)$. If $x\in q(\partial W\times [0,1])$, then $|\tilde
f_1(x)-f(x)|\le(1-t)|f(y)-f(x)|+t|f(x)-f(z_x)|\le2\nabla f(x)$.

Since $f^{-1}(0)\setminus \partial W\cup q((\partial W\cap
f^{-1}(0))\times [0,1])$ is a subcomplex of $f^{-1}(0)$, we have
that $C=\tilde f_1^{-1}(0)=f^{-1}(0)\setminus \partial W\cup
q((\partial W\cap f^{-1}(0))\times [0,1])\cup st(A)$ is a subcomplex
of $N$.

Finally, we modify $\tilde{f}_1$ to obtain the map
$\tilde{f}:N\to[-1,1]$. We push the interior of every $n$ simplex
$\Delta\subset C$ away from $0$ by a small move that is the identity
on $\partial\Delta$. Take the size of the move tending to zero when
$\Delta$ tends to infinity in such a way that the new map $\tilde f$
still has the property that $|f(x)-\tilde f(x)|\le2\nabla f(x)$.
\end{proof}

\begin{Lemma}\label{discrete -> geodesic} Every discrete metric space $Y$ with $\as Y=n\ge 1$
can be isometrically embedded into a geodesic metric space $X$ with
$\as X=n$.
\end{Lemma}

\begin{proof} For every pair of points $x,y\in Y$ we add an interval
$I_{xy}=I_{yx}$ of length $d(x,y)$ and consider the intrinsic metric
on the obtained space.

To prove that $\as X\le n$ we first observe that the collection
$\{I_{xy}\}$ satisfies the inequality $\as \le 1$ uniformly as $x,y$
range over all pairs of points in $Y$. Next observe that for a given
$r$, the set $N_r(Y)\subset X$ is coarsely equivalent to $Y$ and so
has $\as$ $n$. On the other hand, $\{I_{xy}\setminus N_r(Y)\}$ is an
$r$-disjoint collection. So, by the infinite union theorem, $\as
X\le \max\{1,n\}=n.$ Since $Y$ is isometrically embedded in $X$ and
$\as Y=n$, we see that $\as X=n$.
\end{proof}

Next, we need two Lemmas from \cite{DrZ}.

\begin{Lemma}\label{DZ 5.7} Let $f:X\to [-1,1]$ be a continuous map on a geodesic
metric space that is extendible to the Higson corona. Then
$C=f^{-1}(0)$ is an asymptotic separator between $A=f^{-1}(-1)$ and
$B=f^{-1}(1)$.
\end{Lemma}

\begin{proof} This is an immediate consequence of \cite[Lemma 5.7]{DrZ}.
\end{proof}

\begin{Lemma}\label{DZ 5.4}\cite[Lemma 5.4]{DrZ} Suppose that $A$ and $B$ are
asymptotically disjoint subsets of $Y\subset X$ and $C$ is an
asymptotic separator in $X$ with $\as C\le m$. Then there is an
asymptotic separator $\tilde C$ between $A$ and $B$ in $Y$ with $\as
\tilde C\le m$. \qed
\end{Lemma}

A sequence of uniformly bounded locally finite open coverings
$\{\sU_i\}$ of a metric space $X$ is called an {\em anti-\v{C}ech
approximation} \cite{Ro93} for $X$ if $\sU_i\prec\sU_{i+1}$ for all
$i$ and the Lebesgue number $L(\sU_i)$ tends to infinity. Let
$N_i=Nerve(\sU_i)$. Then every anti-\v{C}ech approximation defines a
direct system of locally finite simplicial complexes with refinement
maps as the bonding maps:
\[\begin{CD}
N_1@>q^1_2>> N_2@>q^2_3>>\cdots@>>>N_k@>q^k_{k+1}>>
N_{k+1}@>>>\cdots\hbox{.}
\end{CD}\]

We recall that the projection $p:X\to Nerve(\sU)$ defined by the
partition of unity \[\phi_U(x)=\frac{d(x,X\setminus
U)}{\sum_{V}d(x,X\setminus V)}\] is called the {\em canonical
projection to the nerve}.

An anti-\v{C}ech approximation satisfying the conditions of the
following lemma will be called {\em regular}. Parts (1)-(4) of this
lemma appear as Proposition 1.2 in \cite{Dr05}.

\begin{Lemma}\label{Dr05 Prop1.2} Every proper geodesic metric space $X$ with $\as X\le
n$ admits an anti-\v{C}ech approximation $\{\sU_i,q_{i+1}^i\}$ with
$n$-dimensional locally finite nerves $N_i$ and essentially
surjective projections $p_i:X\to N_i$ such that
\begin{enumerate}
   \item there are bonding maps $p^i_{i+1}:N_i\to N_{i+1}$ with
   $p_{i+1}=p^i_{i+1}\circ p_i$ for all $i$;
   \item $(p^i_{i+1})^{-1}(K)$ is a subcomplex for every subcomplex
   $K\subset N_{i+1}$;
   \item the simplicial maps $q^i_{i+1}:N_i\to N_{i+1}$ are
   simplicial approximations of the $p^i_{i+1}$;
   \item $Lip(p^i_{i+1})<1/2$; and
   \item each $p^i_{i+1}$ is a light map.
\end{enumerate}
\end{Lemma}

\begin{proof} As mentioned above (1)-(4) were proved in
\cite[Proposition 1.2]{Dr05}. To see (5), we observe that each of
the $p_{i+1}^i$ can be chosen to be light, since every map between
$n$-dimensional polyhedra can be approximated by light maps.
\end{proof}

Let $\{p_i:X\to N_n,q^i_{i+1},p^i_{i+1}\}$ be a regular
anti-\v{C}ech approximation of a geodesic metric space $X$. Each map
$q_{i+1}^i$ is simplicial, defined by refinement of covers, and each
nerve $N_i$ is given a uniform geodesic metric $d_i$. So,
$N_i=Nerve(\sU_i)$ and $\sU_i\prec\sU_{i+1}$. Note that $N_i$ is
quasi-isometric to $X$ for every $i$. Let $p_i$ be
$\epsilon_i$-Lipschitz and $m_i$-cobounded, where $\{m_i\}$ is
chosen to be monotone. For $j>i$, we let $p_j^i$ denote
$p_j^{j-1}\circ\cdots\circ p_{i+1}^i$.

Let $N_0$ denote $X$. We define a metric $\bar d$ on the disjoint
union $W=\coprod_{i=0}^{\infty}N_i$ as follows. First, we define
$\bar d$ to be $d_i$ when restricted to $N_i$. Then, $\bar
d(z_i,z_j)=m_j+m_jd_j(p^i_j(z_i),z_j)$ for $z_i\in N_i$, $z_j\in
N_j$ and $i<j$. To prove the triangle inequality, we consider the
projection $\pi: W\to \mathbb{N}_+$, $\pi(N_i)=i$, where
$\mathbb{N}_+$ is given the metric $\rho(i,j)=m_j$ when $i<j$. So,
the metric on $W$ is the sum of the metric on the projection and the
metric on the fiber with the largest index. Since the projection
$p^i_j$ is $1$-Lipschitz, we have the triangle inequality for the
fiber metric when the largest side has a vertex with the maximal
index. The remaining case follows from the inequality $d(x,y)\le
m_kd_k(p_k(x),p_k(y))$.

We are finally in a position to prove that $\asI X\le \as X$ when
$\as X<\infty$.

\begin{proof}[Proof of Theorem \ref{asInd le asdim}]
We use induction on $\as X$. Suppose first that $\as X=0$. Then, by
Theorem 1.1 in \cite{DKU} $\dim\nu X=0$. Thus, $\Ind\nu X=0$. Now,
given two asymptotically disjoint subsets $A$ and $B$ in $X$ we see
that $\emptyset$ is a separator for $A'$ and $B'$ in $\nu X$. This
means that any bounded set is a separator for $A'$ and $B'$. We
conclude that $\asI X\le 0$.

Suppose now that $\as X\le n$ with $n\ge 1$. Given two
asymptotically disjoint subsets $A$ and $B$ in $X$ we have to find
an asymptotic separator $C\subset X$ with $\as C\le n-1$. We may
assume that $X$ is $1$-discrete, otherwise $X$ is coarsely
equivalent to a $1$-discrete space by Proposition \ref{1-discrete}.
By Lemma \ref{discrete -> geodesic}, we may embed $X$ isometrically
into a geodesic metric space $Z$ with $\as Z= n$. Then, we can apply
Lemma \ref{DZ 5.4} to $Z$ and a separator $C'$ in $Z$ to get an
asymptotic separator $C$ in $X$. This means that we may assume that
$X$ is a geodesic space.

Let $A'$ and $B'$ denote the traces in the Higson corona of the
asymptotically disjoint subsets $A$ and $B$ of $X$. Let
$f:\overline{W}\to[-1,1]$ be a continuous function on the Higson
compactification of $W$ such that $f(A)=-1$ and $f(B)=1$. Moreover
we may assume that $f(A_i)=-1$ and $f(B_i)=1$, where $A_i=p_i(A)$
and $B_i=p_i(B)$. Although $A_i$ and $B_i$ are asymptotically
disjoint, they may not be disjoint, so we should remove a compact
$K_i$ from each $N_i$ for all $i$. (This includes the case $A_0$ and
$B_0$, i.e., $A$ and $B$ in $X$.) Formally, we should write
$f:\overline W\setminus \cup_i K_i\to [-1,1]$.

By Proposition \ref{Prop 1}, we may assume that $C_i=f^{-1}(0)\cap
N_i$ is an $(n-1)$-dimensional subcomplex in $N_i$. Indeed, we can
take the star neighborhood $st(C_i)$ along with its regular
neighborhood $V_i$. Take a linear homotopy of $st(C_i)$ to $0$ that
is fixed on the complement to $V_i$. The new map will have the same
extension to the Higson corona and it is unchanged on $A_i$ and
$B_i$. Next, by moves fading to zero, we can push the interiors of
the $n$-simplices from $0$. Then $C_i$ is an $(n-1)$-dimensional
subcomplex of $N_i$ for each $i$.

Let $f'_i=f|_{N_i}$ and $f_i=f'_i\circ p_i$.

We define a separator $V'$ between $A$ and $B$ as the zero-set of a
function $g:X\to [-1,1]$. We construct $g$ by gluing together pieces
of the functions $f_i$. To this end, take a sequence of bounded sets
$R_1\subset V_1\subset R_2\subset V_2\subset\cdots\subset R_i\subset
V_i\subset\cdots$ such that $g$ and $f_i$ agree on $R_i\setminus
V_{i-1}$ where $R_i=p^{-1}(R_i')$ for a subcomplex $R_i'$,
$V_i=p^{-1}(V_i')$, and where $V_i'=st(R_i')$ is the star
neighborhood of $R'_i$. The set $R_{i+1}$ should be viewed as being
much larger than $V_i$.

More formally, we use induction to construct these $R_i$ and
functions $g_i:R_i\to [-1,1]$ such that
\begin{enumerate}
   \item[(a)] $g_{i+1}$ restricted to $R_i$ equals $g_i$; and
   \item[(b)] $g_i$ restricted to $R_i\setminus V_i$ equals $f_i$;
\end{enumerate}
Additionally, we assume that
\begin{enumerate}
   \item[(c)] on $R_i\setminus V_{i-2}$, the function $g_i$ is
obtained by restricting $\tilde g_i\circ p_{i-1}$ for some $\tilde
g_i:N_{i-1}'\to[-1,1]$, where $N_{i-1}'$ is a subcomplex of
$N_{i-1}$ with bounded complement, such that $\tilde{g}_i^{-1}(0)$
is an $(n-1)$-dimensional subcomplex; and
   \item[(d)] $|f-g_i|<8/i$ on $R_i\setminus R_{i-1}$.
\end{enumerate}
From this final condition we see that $g=\cup_ig_i$ will have the
property that $|f-g|\to 0$ at infinity. Thus, $g$ is extendible to
the Higson corona with the same extension as $f$. Then, $g^{-1}(0)$
is an asymptotic separator for $A$ and $B$. Next, we show that
$g^{-1}(0)$ has asymptotic dimension $\le n-1$. To this end, we show
that it admits $1/2^k$-Lipschitz maps to uniform simplicial
complexes of dimension $\le n-1$ for all $k$. Note that
\[p_kg_i^{-1}(0)\subset p_kp_{i-1}^{-1}\tilde
g_{i}^{-1}(0)=(p^k_{i-1})^{-1}\tilde g_i^{-1}(0)\] is an $(n-1)$
dimensional subcomplex of $N_k$ for $k>i$. Denote
$(p^k_{i-1})^{-1}\tilde g_i^{-1}(0)$ by $C_i^k$. Then,
$p_k(g^{-1}(0))=p_k(\cup
g_i^{-1}(0))=p_k(\cup_{i=1}^kg_i^{-1}(0)\cup\cup_{i>k}C^k_i$. Since
the set $p_k(\cup_{i=1}^k g_i^{-1}(0))$ is bounded, collapsing its
star neighborhood to a point gives a uniformly cobounded
$1/2^k$-Lipschitz map of $g^{-1}(0)$ to a uniform
$(n-1)$-dimensional complex.

So, all that remains is to show that the induction goes through. To
this end, we add to the inductive assumptions the following
conditions:
\begin{enumerate}
   \item $|f_{i+1}-f|<1/i$ on $X\setminus R_i$; and
   \item $|\nabla f'_i|,|\nabla f'_{i+1}p^i_{i+1}|<1/i$ on
$N_i\setminus R'_i$.
\end{enumerate}

First, choose $r_i$ so that
\begin{enumerate}
   \item[(1$\empty^\prime$)] $|f_{i+2}-f|<1/i$ on $X\setminus p_i(B_{r_i}(p_i(x_0)))$; and
   \item[(2$\empty^\prime$)] $|\nabla f'_{i+1}|,|\nabla f'_{i+2}p^{i+1}_{i+2}|<1/(i+1)$ on
$N_{i+1}\setminus p^i_{i+1}(B_{r_i}(p_i(x_0)))$.
\end{enumerate}
Then, set $R'_{i+1}=st(p_{i+1}^i(B_{r_i}(p_i(x_0))))$.

Next, define $\tilde g_{i+1}$ by gluing $f_i'$ and
$f_{i+1}'p_{i+1}^i$ along the neighborhood $V_i'\setminus R_i'$. Let
$\phi:N_i\to [0,1]$ be a function so that $\phi(R_i')=0,$ and
$\phi(N_i\setminus V_i')=1$. Not that by the Leibnitz rule and
condition (1), we have $\nabla\overline{f}_{i+1}(z)<2/i$ for
$\overline{f}_{i+1}(z)=\phi(x)f'_{i+1}p^i_{i+1}(z)+(1-\phi(x))f'_i(z)$,
with $z\in N_i\setminus R_i'$. We apply Proposition \ref{Prop 1} to
the function $\overline{f}_{i+1}$ with $A=\mathrm{Int}(V_i'\setminus
R_i')\cap \overline{f}_{i+1}^{-1}(0)$ to obtain a function
$\tilde{f}_{i+1}$ with $\tilde{f}_{i+1}^{-1}(0)$ equal to an
$(n-1)$-dimensional complex. Take $\tilde g_{i+1}$ to be the
restriction of $\tilde{f}_{i+1}$ to a subcomplex $N_i'$ such that
$N_i\setminus N_i'\supset R_i'$. Define
\[g_{i+1}(x)=\left\{
              \begin{array}{ll}
                \tilde{g}_{i+1}(p_i(x)), & \hbox{if $x\in R_{i+1}\setminus R_i$;} \\
                f_i(x), & \hbox{if $x\in R_i$.}
              \end{array}
            \right.
\] For $x\in \partial R_i$, the first function is
$\overline{f}_{i+1}(p_i(x))$ by construction and by Proposition
\ref{Prop 1}. The latter is equal to $f_i(x)$. Hence, the function
$g_{i+1}$ is continuous. Note that \[|f(x)-g){i+1}(x)|\le
|f(x)-f_{i+1}(x)|+|f_{i+1}(x)-g_{i+1}(x)|\le 1/i\] on
$R_{i+1}\setminus V_i$ by the induction hypothesis and the fact that
$g_{i+1}(x)=\tilde{g}_{i+1}p_i(x)=f'_{i+1}p^i_{i+1}p_i(x)=f_{i+1}p_{i+1}(x)=f_{x+1}(x)$
for $x\in R_{i+1}\setminus V_i$.

For $x\in V_i\setminus R_i$, the triangle inequality implies that
$|g_{i+1}(x)-f(x)|\le$\[|\tilde{g}_{i+1}p_i(x)-\overline{f}_{i+1}p_i(x)|
+ |\overline{f}_{i+1}p_i(x)-\overline{f}_{i+1}(z_x)| +
|\overline{f}_{i+1}(z_x)-f'_ip_i(x)| +|f_i(x)-f(x)|\hbox{,}\] where
$d(p_i(x),z_x)\le 1$ and $x_x\in R_i$.

The first summand is no more than $2\nabla\overline{f}_{i+1}p_i\le
4/i$, by (2). The second summand $\le \nabla f_{i+1}p_i(x)\le 2/i$
by (2). The third $\le \nabla f'_ip_i(x)\le 1/i$, by (2). The final
summand is $\le 1/i$ by (1). Thus, $|g_{i+1}(x)-f(x)|\le 8/i$ on
$R_{i+1}\setminus R_i$.
\end{proof}

\section{Hurewicz-type mapping theorem}

In \cite{BD3} the authors prove an asymptotic analog of the Hurewicz
theorem, cf. Theorem \ref{dim Hurewicz Theorem}. In particular, the
following theorem is proved:

\begin{Theorem} \label{BD Hurewicz} Let $f:X\to Y$ be a Lipschitz map from a geodesic
metric space $X$ to a metric space $Y$. Suppose that for every $R>0$
the set family $\{f^{-1}(B_R(y))\}_{y\in Y}$ satisfies the
inequality $\as \le n$ uniformly. Then $\as X\le \as Y+n$.
\end{Theorem}

The proof of the theorem uses mapping cylinders and is too technical
to include in this survey. The main application of this theorem is
to the case where $X$ is a Cayley graph of a finitely generated
group.

Later, Brodskiy, Dydak, Levin and Mitra \cite{BDLM} generalized Bell
and Dranishnikov's result to the following theorem. We follow their
development. First, we give a definition.

Let $f:X\to Y$ be a map of metric spaces. Define \[\as f=\sup\{\as
A\mid A\subset X\hbox{ and }\as(f(A))=0\}\hbox{.}\]

\begin{Theorem} \label{as-Hurewicz} Let $f:X\to Y$ be a large-scale uniform (bornologous
in Roe's terminology) function between metric spaces. Then \[\as
X\le \as Y+\as f\hbox{.}\]
\end{Theorem}

The main idea of the proof is to reformulate the definition of
asymptotic dimension in terms of a double-parameter family that
allows the use of the so-called Kolmogorov Trick. The following is
just a rephrasing of the definition of asymptotic dimension.

\begin{Assertion}Let $X$ be a metric space with $\as X\le n$. Then there is a
function $D_X:\R_+\to\R_+$ such that for each $r>0$ there is there
is a cover $\sU$ of $X$ that can be expressed as
$\bigcup_{i=1}^{n+1}\sU_i$ so that the $\sU_i$ are $r$-disjoint and
$D_X(r)$-bounded.
\end{Assertion}

\begin{defn} The function $D_X$ defined above is called an
\emph{$n$-dimensional control function for $X$}.

Let $k\ge n+1\ge 1.$ An \emph{$(n,k)$-dimensional control function
for $X$} is a function $D_X:\R_+\to\R_+$ such that for any $r>0$
there is a cover $\sU$ of $X$ that can be expressed as
$\bigcup_{i=0}^k\sU_i$ so that the $\sU_i$ are $r$-disjoint and
$D_X(r)$-bounded and such that $\bigcup_{i\in T}\sU_i$ covers $X$
for all subsets $T\subset \{0,1,\dots,k\}$ with $|T|\ge n+1$. Notice
that this condition can be rephrased by saying that each $x\in X$ is
in at least $k-n$ of the $\sU_i$.
\end{defn}

\begin{Lemma} \label{BDLM Lemma} Let $D_X^{n+1}$ be an $n$-dimensional control function
of $X$. Define $\{D_X^i\}_i\ge n+1$ inductively by
\[D_X^{i+1}(r)=D_X^i(3r)+2r\hbox{.}\] Then, each $D^k_X$ is an
$(n,k)$-dimensional control function of $X$ for all $k\ge n+1$.
\end{Lemma}

\begin{proof} We proceed inductively, with the case $k=n+1$ being
trivially true.

Suppose the result to be true for some $k\ge n+1$. Let
$\sU=\bigcup_{i=1}^k \sU_i$ be a $3r$-disjoint, $D_X^k(3r)$-bounded
family so that any $n+1$ of the $\sU_i$ cover $X$. Define $\sU_i'$
to be the $r$-neighborhoods of elements of $\sU_i$ for $i\le k$.
Notice that the elements of $\sU_i'$ are $D^k_X(3r)+2k$-bounded and
$r$-disjoint.

Define $\sU_{k+1}'$ to be the collection of all sets of the form
$\bigcap_{s\in S}A_s\setminus \bigcup_{i\notin S}\sU_i'$, where $S$
is a subset of $\{1,\ldots,k\}$ consisting of exactly $k-n$ elements
and $A_s\in \sU_s$. Observe that each element of $\sU_{k+1}'$ is
contained in some $\sU_j$ so the families $\sU_i$ are
$(D^k_X(3r)+2r)$-bounded.

Next, we must show that the elements of the collection
$\{\sU_i\}_{i=1}^{k+1}$ are $r$-disjoint. Obviously all that needs
to be shown is that the elements of $\sU_{k+1}'$ are $r$-disjoint.
Let $A$ and $B$ be two elements of $\sU_{k+1}'$, say
$A=\bigcap_{s\in S}A_s\setminus \bigcup_{i\notin S}\sU_i'$ and
$B=\bigcap_{t\in T}B_t\setminus \bigcup_{i\notin T}\sU_i'$ with
$S\neq T$. Suppose that $a\in A$ and $b\in B$ with $d(a,b)<r$. Then
there is an $s\in S\setminus T$ such that $a\in A_s$. But, then
there is a $U\in\sU'_s$ containing $b$, a contradiction.

Finally, suppose that $x\in X$ belongs to exactly $k-n$ sets
$\bigcup\sU_i'$, $i\le k$, and let $S=\{i\le k\mid
x\in\bigcup\sU'_i\}$. If $x\notin \bigcup U'_{k+1}$, then $x\in
\bigcup\sU_j'$ for some $j\notin S$, a contradiction. Thus, each $x$
belongs to at least $k+1-n$ elements of $\{\sU_i'\}_{i=1}^{k+1}$.
\end{proof}

Next we prove a product theorem for asymptotic dimension. The
easiest proof (intuitively) involves embedding into uniform
complexes. The product theorem also follows from the asymptotic
Hurewicz Theorem. The following proof comes from \cite{BDLM} and is
a nice illustration of the Kolmogorov Trick.

\begin{Theorem}\label{BDLM product} Let $X$ and $Y$ be metric spaces. Then
\[\as(X\times Y)\le \as X +\as Y\hbox{.}\]
\end{Theorem}

\begin{proof} Put $\as X=m$ and $\as Y=n$ and $k=m+n+1$.
Let $D_X$ be an $(m,k)$-dimension control function for $X$ and $D_Y$
be an $(n,k)$-dimension control function for $Y$. Take $r$-disjoint,
$D_X(r)$-bounded families $\{\sU_i\}_{i=1}^k$ so that any $n+1$
families cover $X$ and $r$-disjoint, $D_Y(r)$-bounded families
$\{\sV_i\}_{i=1}^k$ so that any $m+1$ of the families cover $Y$.
Then the family $\{\sU_i\times \sV_i\}_{i=1}^k$ covers $X\times Y$,
is uniformly bounded and $\sqrt{r}$ disjoint.
\end{proof}
The inequality in this theorem can be strict \cite{BL}. It can be
strict even when one of the factors is the reals $\R$. In
\cite{Dr05} an example of a metric space (uniform simplicial
complex) $X$ is constructed with the properties $\as X=2$ and
$\as(X\times\R)=2$. Thus there is no Morita type theorem for the
asymptotic dimension.

Next we move to the proof of the Hurewicz theorem.

\begin{defn} Let $f:X\to Y$ be a function between metric spaces. We
say that $A\subset X$ is \emph{$(r_X,R_Y)$-bounded} if
$d_Y(f(x),f(x'))\le R_Y$ whenever $d_X(x,x')\le r_X$.
\end{defn}

\begin{defn} Let $f:X\to Y$ be a function of metric spaces. Let
$m\ge 0$. We say that $D_f:\R_+\times \R_+\to \R_+$ is an
\emph{$m$-dimensional control function of $f$} if for any $r_X>0$
and $R_Y>0$ and any $A\subset X$ with $\diam(f(A))\le R_Y$, $A$ can
be expressed as the union of $m+1$ sets whose $r_X$-components are
$D_f(r_X,R_Y)$-bounded.
\end{defn}

\begin{prop} \label{BDLM 4.5} Suppose that $f:X\to Y$ is a function between metric
spaces and that $m\ge 0$. If $\as f\le m$ then $f$ has an
$m$-dimensional control function $D_f$.
\end{prop}

\begin{proof} Fix non-negative $r_X$ and $R_Y$. Suppose that for
each $n$ there is a $y_n\in Y$ such that $A_n=f^{-1}(B(y_n,R_Y))$
cannot be expressed as a union of $m+1$ sets whose $r_X$-components
are $n$-bounded. Then, the set $C=\bigcup_{n=1}^\infty B(y_n,R_Y)$
cannot be bounded. If $C$ were bounded, then $\as f^{-1}(C)\le m$.
By passing to a subsequence, we may arrange $y_n\to\infty$ and $\as
C=0$, a contradiction.
\end{proof}

\begin{defn} Let $k\ge M+1\ge 1$. An \emph{$(m,k)$-dimensional control
function of $f:X\to Y$} is a function $D_f:\R_+\times\R_+\to \R_+$
such that for all $r_X>0$ and all $R_Y>0$ any $(\infty,R_Y)$-bounded
subset $A\subset X$ can be expressed as the union of $k$ sets
$\{A_i\}_{i=1}^k$ whose $r_X$-components are $D_f(r_X,R_Y)$-bounded
and so that any $x\in A$ belongs to at least $k-m$ of the $A_i$.
\end{defn}

\begin{prop}\label{BDLM 4.7} Let $f:X\to Y$ be a function between
metric spaces and $m\ge 0$. Suppose that there is an $m$-dimensional
control function $D^{m+1}_f$ of $f$. Define $D^{k}_f:\R_+\times
\R_+\to\R_+$ inductively by
$D^{k}_f(r_X,R_Y)=D^{k-1}_f(3r_X,R_Y)+2r_X$ for each $k\ge m$. Then
each $D^k_f$ is an $(m,k)$-dimensional control function of $f$.
\end{prop}

The proof is similar to that of Lemma \ref{BDLM Lemma}.

\begin{Theorem} \label{BDLM 4.9} Let $k=m+n+1$ and suppose $f:X\to Y$ is a
large-scale uniform function of metric spaces with $\as Y\le n$. If
there is an $(m,k)$-dimensional control function $D_f$ of $f$ then
$\as X\le m+n$.
\end{Theorem}

The proof is a bit technical, so we give the idea; the interested
reader is referred to the original article \cite{BDLM}.

\begin{proof}[Sketch of proof.] Given $r$ we want to find a constant $D_X(r)$ so that
$X$ can be written as the union $X=\bigcup_{j=1}^k D^j$ with each
$r$-component $D_X(r)$-bounded. This clearly implies the conclusion
of the Theorem.

We want to apply the Kolmogorov Trick, so we take a cover of $Y$ by
sets $A_i$, ($i=1,\dots,k$) with conditions on boundedness of
components. For each $i$ write $A_i$ as a union of $k$ sets
$\{U^j_i\}$ with certain boundedness conditions on components in
such a way that every point of $y$ belongs to at least $m$ sets. For
each $i$ cover the set $f^{-1}(A_i)$ by $k$ sets $\{B^j_i\}_{j=1}^k$
with conditions on boundedness of components so that every point in
$f^{-1}(A_i)$ is in at least $m$ sets. Then, apply the Kolmogorov
Trick to $D^j=\bigcup_i\left(B^j_i\cap f^{-1}(U^j_i)\right)$.
\end{proof}

\begin{Corollary} Suppose $f:X\to Y$ is a function of metric spaces
and $m\ge 0$. Then, $\as f\le m$ if and only if $f$ has an
$m$-dimensional control function.
\end{Corollary}

\begin{proof} The ``only if" part is Proposition \ref{BDLM 4.5}. The
``if" part follows from Theorem \ref{BDLM 4.9} with $n=0$.
\end{proof}

\begin{Theorem} Suppose that $f$ is a large-scale uniform function
$f:X\to Y$ between metric spaces. Then \[\as X\le \as Y+\as
f\hbox{.}\]
\end{Theorem}

\begin{proof} Put $n=\as Y$ and $m=\as f$. Then by Theorem \ref{BDLM
4.9} it suffices to show that there is an $(m,k)$-dimensional
control function $D_f$ for $f$ with $k=m+n+1$. By the previous
corollary, $f$ has an $m$-dimensional control function. Finally, by
Proposition \ref{BDLM 4.7} there is an $(m,k)$-dimensional control
function, and the proof is complete.
\end{proof}

\section{Coarse Embedding}

The goal of this section is to see that a metric space with finite
asymptotic dimension admits a coarse embedding into Hilbert space.
This result is of particular interest in connection with the Novikov
higher signature conjecture. Guoliang Yu showed in \cite{Yu1} that
finitely generated groups that admit a coarse embedding into Hilbert
space satisfy this conjecture.

In \cite{Yu2}, Yu defined a property called ``Property A," which can
be thought of as a generalization of amenability for discrete spaces
with bounded geometry. For groups, this is equivalent to the
exactness of the reduced $C^*$-algebra by a result of Ozawa
\cite{Oz}, so it has come to be known as exactness of the group.
This property is implied by finite asymptotic dimension and is
enough to imply a coarse embedding into Hilbert space.

\begin{Theorem} Let $X$ be a metric space with finite asymptotic
dimension. Then $X$ admits a coarse embedding into Hilbert space.
\end{Theorem}

\begin{proof} We follow Roe's proof in \cite{Ro03}.
By the fifth characterization of asymptotic dimension in Theorem
\ref{equiv asdim}, for every $k>0$ we can find uniformly cobounded
$2^{-k}$-Lipschitz maps $\phi_k:X\to L_k$ where $L$ is a finite
dimensional simplicial complex with the metric induced from its
inclusion in Hilbert space. Let $x_0$ be some fixed basepoint in
$X$. Then define $\Phi:X\to \bigoplus_{k=1}^\infty \sH$ by
$\Phi(x)=\{\phi_k(x)-\phi_k(x_0)\}.$

Then, since $\phi_k$ is $2^{-k}$-Lipschitz, we see that
\[\|\Phi(x)-\Phi(x')\|^2=\sum_{k=1}^\infty\|\phi_k(x)-\phi_k(x')\|^2\le
\|x-x'\|^2\sum_{k=1}^\infty 2^{-2k}.\]

On the other hand, since each $\phi_k$ is uniformly cobounded, this
means there is a $R_k$ so that $\diam(\phi_K^{-1}(\sigma))\le R_k$
for each simplex $\sigma\in L_k.$ So if $\|x-x'\|^2>R_\ell$ for some
$\ell$ then $\phi_k(x)$ and $\phi_k(x')$ are orthogonal unit vectors
for $k\le \ell.$ Thus, $\|\phi_k(x)-\phi_k(x')\|^2=2$ for all $k\le
\ell$ and so $\|\Phi(x)-\Phi(x')\|^2\ge 2\ell.$
\end{proof}

Note that there cannot be an asymptotic analog of the
N\"obeling-Pontryagin Theorem (Theorem \ref{NP}). For example, the
binary tree $T$ has $\as T=1$, but it does not admit a coarse
embedding in Euclidean space of any dimension. The obstacle is
volume growth. Nevertheless there is the following analog of Theorem
\ref{Thm17} \cite{Dr3}.
\begin{Theorem} Every proper metric space $X$ with $\as X\le n$ admits a
coarse embedding in $n+1$ locally finite trees.
\end{Theorem}

In \cite{DrZ} a universal metric space for the class of metric
spaces of bounded geometry and asymptotic dimension $\le n$ is
constructed. It is not an asymptotic analog of the Menger space
$\mu^n$ since it does not have bounded geometry. Moreover, it is
proven in \cite{DrZ} that there is no proper universal space for
asymptotic dimension $n$. Macro-scale analogs of N\"obeling spaces
have been constructed that are universal for asymptotic dimension
and coarse embeddings \cite{BN}.

\section{Hyperbolic spaces}

Let $X$ be a metric space. For $x$, $y$, $z\in X$ we define the {\em
Gromov product}
$$(x|y)_z:=\frac{1}{2}(|zx|+|zy|-|xy|).$$
Let $\delta\ge 0$. A triple $(a_1,a_2,a_3)\in\R^3$ is called a {\em
$\delta$-triple}, if $a_{\mu}\ge\min\{a_{\mu+1},a_{\mu+2}\}-\delta$
for $\mu =1,2,3$, where the indices are taken modulo 3.

The space $X$ is called {\em hyperbolic} if there is $\delta>0$ such
that for every $o$, $x$, $y$, $z\in X$ the triple
$((x|y)_o,(y|z)_o,(x|z)_o)$ is a $\delta$-triple.

Note that if $X$ satisfies the $\delta$-inequality for one
individual base point $o\in X$, then it satisfies the
$2\delta$-inequality for any other base point $o' \in X$, see, for
example \cite{Gr87}. Thus, to check hyperbolicity, one has to check
this inequality only at one point.

Let $X$ be a hyperbolic space and $o\in X$ be a base point. A
sequence of points $\{x_i\}\subset X$ {\em converges to infinity,}
if $\lim_{i,j\to\infty}(x_i|x_j)_o=\infty.$ Two sequences $\{x_i\}$,
$\{x_i'\}$ that converge to infinity are {\em equivalent} if
$\lim_{i\to\infty}(x_i|x_i')_o=\infty.$ Using the
$\delta$-inequality, one easily sees that this defines an
equivalence relation for sequences in $X$ converging to infinity.
The {\em boundary at infinity} $\partial_{\infty}X$ of $X$ is
defined as the set of equivalence classes of sequences converging to
infinity.

A hyperbolic space $Y$ is said to be {\em visual}, if for some base
point $o\in Y$ there is a positive constant $D$ such that for every
$y\in Y$ there is $\xi\in\di Y$ with $|oy|\le(y|\xi)_o+D$ (one
easily sees that this property is independent of the choice of $o$).
Here $(y|\xi)_o=\inf\liminf_{i\to\infty}(y|x_i)_o,$ where the
infimum is taken over all sequences $\{x_i\}\in\xi$. For hyperbolic
geodesic spaces this property is a rough version of the property
that every segment $oy\subset Y$ can be extended to a geodesic ray
beyond the end point $y$.

Most of our interest is in geodesic metric spaces. A geodesic metric
space is hyperbolic if there is $\delta>0$ such that every geodesic
triangle is $\delta$-thin \cite{BH}, which means that every side of
the triangle is contained in a $\delta$-neighborhood of the other
two.

The following is straightforward.
\begin{Proposition}
If a geodesic metric space is quasi-isometric to hyperbolic space
then it is hyperbolic.
\end{Proposition}

In the second part of this paper, we prove that the asymptotic
dimension of a finitely generated hyperbolic group is finite. In
fact, this applies to more general hyperbolic spaces, see for
example \cite{BellFujiwara,BS,Ro05}. We note that in the case of
hyperbolic groups $\as\Gamma=\dim\partial_{\infty}\Gamma+1$
\cite{BL},\cite{Bu}.

The fundamental group $\pi_1(X,x_0)$ of a metric space $X$ is
generated by a set of free loops $S$ if it is generated by the set
of based loops of the form $\phi=pf\bar p$, $f\in S$,  where $p$ is
a path from $x_0$ to $f(0)=f(1)$.

The fundamental group $\pi_1(X)$ of a metric space $X$ is {\it
uniformly generated} if there is $L>0$ such that $\pi_1(X)$ is
generated by free loops of length $\le L$.

The following proposition can be extracted from \cite{FW}.
\begin{Proposition}
The fundamental group of a hyperbolic space $X$ is uniformly
generated.
\end{Proposition}
In the case of hyperbolic spaces the embedding theorem of Section 8
can be improved to the following \cite{BDS}
\begin{Theorem} Every visual hyperbolic space
$X$ admits a quasi-isometric embedding into the product of $n+1$
copies of the binary metric tree where $n=\dim\di X$ is the
topological dimension of the boundary at infinity.
\end{Theorem}

\section{Spaces of dimension 0 and 1}

Recall that spaces with asymptotic dimension $0$ are those spaces
that can be presented as a union of uniformly bounded, $r$-disjoint
sets for each (large) $r$. Notice that the collection of
asymptotically zero-dimensional metric spaces includes all compacta,
but is clearly not limited to such things. On the other hand, if $X$
is assumed to be geodesic, then $\as X=0$ implies that $X$ is
compact. For asymptotically $1$-dimensional spaces we have the
following result.

\begin{Theorem} \label{asdim1}Let $X$ be a geodesic metric space with $\as X=1$
whose fundamental group is uniformly generated. Then $X$ is
quasi-isometric to an infinite tree.
\end{Theorem}

\begin{proof}
Since $\as X=1$, there is a $1/2L$-Lipschitz map $p:X\to K$ to a
uniform 1-dimensional simplicial complex which is a quasi-isometry.
Thus, $p$ sends every loop of length $\le L$ to a null-homotopic
loop. Therefore, $p_*:\pi_1(X)\to\pi_1(K)$ is the zero homomorphism.
Since every geodesic space is path connected and locally path
connected, by the Lifting Criterion there is a lift $\tilde
p:X\to\tilde K$ of $p$ to the universal cover $u:\tilde K\to K$.
Since $u$ is a local isometry, $\tilde p$ is locally Lipschitz.
Since $X$ is geodesic $\tilde p$ is globally Lipschitz. Clearly, it
is a quasi-isometry onto the image $T=\tilde p(X)$. Note that $T$ is
a tree as a connected subcomplex of a tree $\tilde K$.
\end{proof}

This result first appeared as \cite[Theorem 0.1]{FW}. There, the
proof appeals to Manning's bottle-neck property, \cite{Manning}. The
main application of this result is the following
\begin{theorem}[\cite{FW}] Let $S$ be a compact oriented surface
of genus $g\ge 2$ and with one boundary component. Let $C(S)$ be the
curve graph of $S$. Then $\as C(S)>1$.
\end{theorem}
\begin{proof}
In this case $C(S)$ is one-ended by a result of Schleimer
\cite{Sch}. Hence it cannot be quasi-isometric to a tree. Since
$C(S)$ is hyperbolic \cite{MazMin}, we obtain a contradiction with
Theorems 48 and Proposition 46.
\end{proof}
Here we recall that curve graph of a surface $S$ is the graph whose
vertices are isotopy classes of essential, nonperipheral, simple
closed curves in $S$, with two distinct vertices joined by an edge
if the corresponding classes can be represented by disjoint curves.
Bell and Fujiwara proved that $\as C(S)<\infty$ for all surfaces $S$
\cite{BellFujiwara}.

\section{Linear control}

Let $X$ be a metric space with $\as X\le n$. If there is a $C>0$ so
that for every $D$, there is a cover $\sU=\sU_0\cup\cdots\cup\sU_n$
of $X$ by $D$-disjoint sets with $\mesh(\sU)<CD$, then we say that
$\as X\le n$ with {\it linear control}. This property was defined in
\cite{DrZ}, where it was called the \textit{Higson property}. This
is related to the \textit{Assouad-Nagata dimension}, which is
defined as follows.

\begin{defn} For a metric space $X$, the
{\it Assouad-Nagata dimension} $AN-\dim X$ is the infimum of all
integers $n$ such that there is a $C>0$ so that for any $D>0$, $X$
can be covered by a $CD$-bounded cover with $D$-multiplicity $\le
n+1$ \cite{A}.
\end{defn}
Using the Assouad-Nagata dimension U. Lang and T. Schlichenmaier
gave the following refinement of Theorem 44 \cite{La}:

\begin{Theorem} If for a metric space AN-$\dim (X,d)\le n$ then
for sufficiently small $\epsilon$, $(X,d^{\epsilon})$ admits a
bi-Lipschitz embedding in the product of $n+1$ locally finite trees.
\end{Theorem}
The Assouad-Nagata dimension is a way of simultaneously considering
dimension at all scales. In \cite{BDHM} the Assouad-Nagata dimension
was characterized in terms of extension of Lipschitz maps to the
unit $n$-sphere $S^n$. When it is applied to discrete spaces like
finitely generated groups the small scales are not important and it
defines a quasi-isometry invariant. Many of the results contained in
this survey have corresponding results in the theory of
Assouad-Nagata dimension, (eg. in \cite{BDLM} a Hurewicz-type
theorem for Assouad-Nagata dimension is proved).

Although a thorough discussion of Assouad-Nagata dimension is beyond
the scope of this survey, we do want to draw attention to some
curious features of linear control.

\begin{prop}\cite[Proposition 4.1]{DrZ} Every proper metric space
$X$ is coarsely equivalent to a proper space $Y$ with $\as Y\le n$
with linear control.
\end{prop}

The proof of this proposition in \cite{DrZ} uses universal spaces
for asymptotic dimension. Each metric space admits a coarsely
uniform embedding into such a universal space, and the universal
space has asymptotic dimension $\le n$ with linear control. Since
linear control passes to subsets, the (coarsely) embedded copy of
$X$ has $\as\le n$ with linear control. An alternative proof is
given in \cite{BDHM} where it was shown that $Y$ can be taken to be
hyperbolic. In case of groups it was shown in \cite{BDLM} that there
is a proper left-invariant metric for which $\as G=\textit{AN-}\dim
G$.

This is an illustration of a difference between quasi-isometry and
coarse equivalence. Whereas quasi-isometry will preserve the
property of linear control, coarse equivalence will not.

In particular, the Morita type theorem holds for Assouad-Nagata
dimension of cocompact spaces \cite{DrSm2}:
\begin{theorem}
$AN{-}\dim(X\times\R)=AN{-}\dim X+1$.
\end{theorem}

In \cite{BDL} it was shown that the lamplighter group
$G=\Z_2\wr\Z^2$ has $\as G=2$ and $AN{-}\dim G=\infty$.

In \cite{DHi} it was proven that the dimension of asymptotic cone of
a metric space does not exceed its Assouad-Nagata dimension.:
\begin{theorem}
$\dim(cone_{\omega}X)\le AN{-}\dim X$ for all ultrafilters
$\omega\in\beta\mathbf{N}$.
\end{theorem}
We recall the definition of the asymptotic cone $cone_{\omega}(X)$
of a metric space with base point $x_0\in X$ with respect to a
non-principal ultrafilter $\omega$ on $\mathbf{N}$
\cite{Gr93},\cite{Ro03}. On the sequences of points $\{x_n\}$ with
$\|x_n\|\le Cn$ for some $C$, we define an equivalence relation
$$\{x_n\}\sim\{y_n\}\Leftrightarrow \lim_{\omega}d(x_n,y_n)/n=0.$$
We denote by $[\{x_n\}]$ the equivalence class of $\{x_n\}$. The
space $cone_{\omega}(X)$ is the set of equivalence classes
$[\{x_n\}]$ with the metric
$d_{\omega}([\{x_n\}],[\{y_n\}])=\lim_{\omega}d(x_n,y_n)/n$. We note
that the space $cone_{\omega}(X)$ does not depend on the choice of
the base point.

\section{Dimension of general coarse structures}

So far we have only considered the asymptotic dimension of metric
spaces. John Roe has shown that what is really important for this
large-scale version of dimension is the so-called coarse structure
of the space. This should be thought of as the analog of the
situation one encounters beginning a study of topology. The first
examples one sees in topology are metric spaces, but one soon
realizes that the abstract notion of a topology is really what makes
the theory work. This section follows the development given in
\cite{Ro03}. An alternative approach to the coarse structures is
given in \cite{DH} and some of the basic properties of asymptotic
dimension in the coarse sense are developed in \cite{Grave}.

To begin, we define the abstract notion of a coarse structure. Let
$X$ be a set. If $E\subset X\times X$, then the \emph{inverse} of
$E$, denoted $E^{-1}$ is the set $\set{(x,x')\mid (x',x)\in E}$. If
$E'$ and $E''$ are subsets of $X\times X$, then the product is
denoted $E'\circ E''$ and is defined to be
\[E'\circ E''=\set{(x',x'')\mid \exists x\in X\hbox{ such that }
\exists(x',x)\in E'\hbox{ and }\exists (x,x'')\in E''}\hbox{.}\]

If $E\subset X\times X$ and $K\subset X$, define $E[K]=\{x'\mid
\exists\, x\in K, (x',x)\in E\}$. When $K={x}$ we use the notation
$E_x=E[\{x\}]$ and $E^x=E^{-1}[\{x\}]$.

\begin{defn} A {\em coarse structure} on a set $X$ is a collection
$\mathcal{E}$ of sets (called {\em controlled sets} or {\em
entourages} for the coarse structure) that contains the diagonal and
is closed under the formation of subsets, inverses, products and
finite unions. A set equipped with a coarse structure is called a
{\em coarse space.}
\end{defn}
A subset $D\subset X$ of a coarse space is \emph{bounded} if
$D\times D$ is controlled. A family of sets $\{D_i\}$ is
\emph{uniformly bounded} if $\cup_i(D_i\times D_i)$ is controlled.

Recall the following characterization of asymptotic dimension for
metric spaces:
\begin{enumerate}
   \item $\as X=0$ on $r$-scale if and only if there
   is a cover of $X$ by uniformly bounded, $r$-disjoint sets.
   \item $\as X\le n$ if and only if, for every $r<\infty$,
   $X$ can be written as a
   union of at most $n+1$ sets that are $0$-dimensional on $r$-scale.
\end{enumerate}

To translate this to the coarse category, we need a notion of
disjointness.

\begin{defn} Let $X$ be a coarse space and let $U$ be a controlled
set. We say that $D\subset X$ is $U$-{\em disconnected} if it can be
written as a disjoint union $D=\bigsqcup_{i=1}^\infty D_i$ such that
\begin{enumerate}
   \item $\{D_i\}$ is uniformly bounded, i.e., $\cup (D_i\times D_i)=W$ is controlled,
   (in this case, we call the family $\{D_i\}$ $W$-{\em bounded}) and
   \item when $i\neq j$, $D_i\times D_j$ is disjoint from $U$ (in this case, we call the family
   $\{D_i\}$ $U$-{\em disjoint}).
\end{enumerate}
\end{defn}

\begin{defn} Let $(X,\mathcal{E})$ be a coarse space. Then:
\begin{enumerate}
   \item $\as X=0$ if it is $U$-disconnected for every controlled set
$U$.
   \item $\as X\le n$ if for every controlled set $U$, $X$ can be
written as the union of at most $n+1$ $U$-disconnected subsets.
\end{enumerate}
\end{defn}

\begin{example} Let $(X,d)$ be a metric space and let $\mathcal{E}$
be the collection of all $E\subset X\times X$ for which
$\sup\{d(x,x')\mid (x,x')\in E\}$ is finite. Then, $\mathcal{E}$ is
a coarse structure called the {\em bounded coarse structure
associated to $(X,d)$}.
\end{example}

It is straightforward to prove the following proposition.

\begin{prop} Let $(X,d)$ be a metric space with bounded coarse structure
$\mathcal{E}$. Then $\as (X,d)=\as (X,\mathcal{E})$.\qed
\end{prop}

A coarse structure on a topological space is {\em consistent with
the topology} if the bounded sets for this structure are exactly
those that are relatively compact. Suppose $\mathcal{E}$ is a coarse
structure that is consistent with the topology on a locally compact
space $X$. We say that $f:X\to\mathbf{C}$ is a {\em Higson
function}, denoted $f\in C_h(X,\mathcal{E})$ if for every
$E\in\mathcal{E}$ and every $\epsilon>0$, there is a compact set
$K\subset X$ such that $|f(x)-f(y)|<\epsilon$ whenever $(x,y)\in
E\setminus(K\times K)$. Then by the Gelfand-Naimark theorem there is
a compactification $h_{\mathcal{E}}X$ of $X$ called the {\em Higson
compactification}. The {\em Higson corona} is defined by
$\nu_{\mathcal{E}}X=h_{\mathcal{E}}X\setminus X$.

We consider a proper metric space $(X, d)$ with basepoint $x_0$ and
define $\norm{x} = d(x, x_0)$.  

\begin{definition} We define the {\it sublinear coarse structure}, denoted
$\mathcal{E}_{L}$,

on $X$ as follows:
$$ \mathcal{E}_{L} = \{ E \subset X\times X :
\lim_{x \rightarrow \infty} \frac{\sup_{y \in E_x} d(y,x)
}{\norm{x}} = 0 = \lim_{x \rightarrow \infty} \frac{\sup_{y \in E^x}
d(x,y) }{\norm{x}} \}. $$ By the statement $\lim_{x \rightarrow
\infty} \frac{\sup_{y \in E_x} d(y,x) }{\norm{x}} = 0$, we mean that
for each $\epsilon > 0$, there is a compact subset $K$ of $X$
containing $x_0$ (equivalently, an $r\geq 0$) such that
$$ \frac{\sup_{y \in E_x} d(y,x) }{\norm{x}} \leq \epsilon $$
for all $x \notin K$ (respectively, for all $x$ with $\norm{x} >
r$).  It would perhaps be better to think of this as $\lim_{\norm{x}
\rightarrow \infty}$. We leave to the reader to check that
$\mathcal{E}_{L}$ is indeed a coarse structure and it does not
depend on the choice of basepoint. The Higson corona for the
sublinear coarse structure on $X$ will be denoted by $\nu_LX$.
\end{definition}

The sublinear coarse structure is useful for the Assouad-Nagata
dimension \cite{DrSm2}
\begin{theorem}
For a cocompact connected proper metric space,
$$AN{-}\dim{X}=\dim\nu_LX$$ provided $AN{-}\dim{X}<\infty$.
\end{theorem}
A metric space $X$ is {\em cocompact} if there is a compact set
$C\subset X$ such that $Iso(X)(C)=X$ where $Iso(X)$ is the group of
isometries of $X$.

\

{\bf II. ASYMPTOTIC DIMENSION OF GROUPS}

\

\section{Metrics on groups}
A \textit{norm} on a group $G$ is a unary operation $\|\cdot\|$
satisfying
\begin{enumerate}
   \item $\|g\|=0$ if and only if $g=e$;
   \item $\|g\|=\|g^{-1}\|$ for all $g\in G$; and
   \item $\|gh\|\le \|g\|+\|h\|.$
\end{enumerate}
Let $\Gamma$ be a finitely generated group. To any finite generating
set $S=S^{-1},$ one can assign a norm $\|\cdot\|_S$ defined by
setting $\|g\|_S$ equal to the length of the shortest $S$-word
presenting the element $g.$

One can now define the \textit{left-invariant word metric associated
to $S$} \index{word metric} by $d_S(g,h)=\|g^{-1}h\|_S$. When the
generating set is understood we write $d(g,h)$ for $d_S(g,h)$ and
$\|\cdot\|$ for $\|\cdot\|_S$. This metric is left-invariant, i.e.,
the action of $G$ on itself by left multiplication is an isometry:
$d(ag,ah)=\|g^{-1}a^{-1}ah\|=\|g^{-1}h\|=d(g,h).$  The word metric
turns the group $\Gamma$ into a discrete metric space with bounded
geometry. Since closed balls are finite -- and hence compact -- a
finitely generated group in the word metric is proper.

We define the asymptotic dimension of a finitely generated group
$\Gamma$ by $\as\Gamma=\as(\Gamma,d_S),$ where $S$ is any finite,
symmetric generating set for $\Gamma.$

\begin{Corollary} \label{inv} Let $\Gamma$ be a finitely generated group.  Then
$\as\Gamma$ is an invariant of the choice of generating set, i.e.,
it is a group property.
\end{Corollary}

\begin{proof} Let $S$ and $S'$ be finite generating sets for
$\Gamma.$  We have to show that $(\Gamma,d_S)$ and $(\Gamma,d_{S'})$
are coarsely equivalent.  In fact, they are Lipschitz equivalent, as
we now show.

Let $\lambda_1=\max\{\|s\|_{S'}\mid s\in S\}$ and
$\lambda_2=\max\{\|s'\|_S\mid s'\in S'\}.$  It follows that
$\lambda_2^{-1}\|\gamma\|_{S'}\le\|\gamma\|_S\le
\lambda_2\|\gamma\|_{S'}.$  Take
$\lambda=\max\{\lambda_1,\lambda_2\}.$ Then $\lambda^{-1}
d_{S'}(g,h)\le d_S(g,h)\le \lambda d_{S'}(g,h).$
\end{proof}

%

This proof shows that from the large-scale point of view two word
metrics on a finitely generated group $\Gamma$ are
indistinguishable.  

Alternatively, we could define $\as\Gamma$ to be $\as
\left(Cay(\Gamma,S)\right)$ where $Cay(\Gamma,S)$ denotes the Cayley
graph of $\Gamma$ with respect to the generating set $S.$ It is easy
to see that $\Gamma$ with the word metric associated to $S$ and the
Cayley graph $Cay(\Gamma,S)$ with its edge-length metric are
quasi-isometric.

\subsection{Coarse Equivalence on Groups}

For finitely generated groups, the primary notion of large-scale
equivalence is that of quasi-isometry. Since finitely generated
groups are quasi-isometric to geodesic metric spaces (Cayley
graphs), a coarse equivalence between them is a quasi-isometry. A
coarse equivalence is practically useful when dealing with
countable but not finitely generated groups. We note that an
alternative (equivalent) approach was implemented by means of
extension of the notion of quasi-isometry to all countable groups
in \cite{Sh}.

A fundamental result in geometric group theory is the following. 

\begin{theorem}[\v{S}varc-Milnor Lemma]\label{Svarc-Milnor}
\index{Svarc-Milnor Lemma@\v{S}varc-Milnor Lemma} Let $X$ be a
proper geodesic metric space and $\Gamma$ a group acting properly
cocompactly by isometries on $X$. Then $\Gamma$ is finitely
generated and (in the word metric) $\Gamma$ and $X$ are
quasi-isometric. In particular, they are coarsely equivalent.
\end{theorem}

\begin{cor} The following are easy consequences of the
\v{S}varc-Milnor Lemma:
\begin{enumerate}
   \item Let $\Gamma$ be a finitely generated group and
   let $\Gamma'\subset\Gamma$ be a finite index subgroup. Then
   $\Gamma$ and $\Gamma'$ are quasi-isometric.
   \item Let \[\begin{CD}1@>>>K@>>>\Gamma@>\phi>>H@>>>1\end{CD}\]
   be an exact sequence with
   $K$ finite and $H$ finitely generated. Then $\Gamma$ is finitely
   generated and quasi-isometric to $H$.
\end{enumerate}\qed
\end{cor}



Since asymptotic dimension is an invariant of quasi-isometry, we
immediately obtain:

\begin{cor} Let $\Gamma$ be a finitely generated group.
\begin{enumerate}
   \item  Let $\Gamma'\subset\Gamma$ be a finite index subgroup. Then
   $\as\Gamma=\as\Gamma'$.
   \item Let \[\begin{CD}1@>>> K@>>>\Gamma@>>>H@>>>1\end{CD}\] be an exact sequence with
   $K$ finite and $H$ finitely generated. Then $\as\Gamma=\as H$
   (cf. Theorem \ref{extension}).
\end{enumerate}\qed
\end{cor}

Two groups are said to be {\em commensurable} if they have
isomorphic finite-index subgroups. The previous corollary
immediately implies the following result.

\begin{cor} Let $\Gamma$ and $\Gamma'$ be commensurable with $\Gamma$ finitely
generated. Then $\as\Gamma=\as\Gamma'$.\qed
\end{cor}


Another application of the \v{S}varc-Milnor Theorem is the
following, cf. \cite{Gr93}.

\begin{cor}\label{cor57}  Let $M$ be a compact Riemannian manifold with universal
cover $\tilde{M}$ and (finitely generated) fundamental group $\pi$.
Then $\as \tilde{M}=\as\pi$.
\end{cor}

\subsection{Asymptotic Dimension of General Groups}

An immediate problem that one encounters when dealing exclusively
with finitely generated groups is that a finitely generated group
can have a subgroup that is not finitely generated. Although it is
tempting to define the asymptotic dimension of the subgroup to be
the asymptotic dimension of the metric subspace of the finitely
generated group, the asymptotic dimension of the subgroup should be
defined in terms of its abstract group structure, not in terms of a
particular homomorphic embedding into a larger, finitely generated
group.

It is for this reason that we make the following convention. When
dealing with a countable, possibly non-finitely generated group $G$
we define the asymptotic dimension $\as G$ to be the asymptotic
dimension of the group endowed with some left-invariant, proper
metric. It remains to show that the choice of left-invariant proper
metric does not affect the asymptotic dimension. This is guaranteed
in view of the following \cite{Sh}, \cite[Proposition 1]{Smi}
\begin{Proposition} \label{proper coarse}
Let $\Gamma$ be a countable group. Then any two left-invariant
proper metrics on $\Gamma$ are coarsely equivalent.
\end{Proposition}

\begin{cor} \label{subgroup cor} Let $\Gamma$ be a finitely generated group and
let $\Gamma'\subset \Gamma$. Then, $\as \Gamma'\le\as\Gamma$.
\end{cor}

\begin{proof} Since all left-invariant proper metrics on $\Gamma'$
are coarsely equivalent, we need only consider the asymptotic
dimension of $\Gamma'$ as a subspace of $\Gamma$, where $\Gamma$ is
endowed with the left-invariant word metric associated to a finite
generating set. The fact that $\as \Gamma'\le\as\Gamma$ as a
subspace is an easy consequence of the definition.
\end{proof}
Moreover, the following holds \cite{DrSm1}
\begin{Theorem} Let $G$ be a countable group. Then $\as
G=\sup\{\as F\mid F\subset G \hbox{ is finitely generated}\}$.
\end{Theorem}

This leads to the following
\begin{defn} Let G be a (possibly uncountable) group
\[\as G=\sup\{\as F\mid F\subset G\hbox{ is finitely generated}\}\hbox{.}\]
\end{defn}

Notice that with this definition, it is immediate that $\as H\le \as
G$ for any subgroup $H\subset G$.

There is a remark in \cite{DrSm1} that the definition of asymptotic
dimension of arbitrary groups given above coincides with the
asymptotic dimension of the coarse space $(G,\mathcal{E})$ where
$E\in\mathcal{E}$ if and only if $\{x^{-1}y\mid (x,y)\in E\}$ is
finite.

It should also be noted that for an uncountable group, the
asymptotic dimension as defined here does not necessarily agree with
the asymptotic dimension the group of equipped with a left-invariant
proper metric.

Taking groups with proper left-invariant metrics can give rise to
interesting examples, as pointed out by \cite{Smi}.

\begin{example} Endow $\mathbb{Q}$ with a left-invariant proper
metric. Then $\as\mathbb{Q}=0.$ On the other hand, as a metric
subspace of $\mathbb{R},$ $\mathbb{Q}$ is coarsely equivalent to
$\mathbb{R}$ so $\as\mathbb{Q}=1.$ The problem is that the metric
$\mathbb{Q}$ inherits from $\mathbb{R}$ is not a proper metric.
\end{example}

\subsection{Groups with asdim 0 and 1}

In the first section we proved stated that geodesic spaces with
asymptotic dimension $0$ are compact. This easily implies the
following:

\begin{prop} Let $\Gamma$ be a finitely generated group. Then
$\as\Gamma=0$ if and only if $\Gamma$ is finite.
\end{prop}

The following theorem was proven independently by several authors
\cite{Ge} \cite{JS}.
\begin{Theorem}
Every finitely presented group $\Gamma$ with $\as\Gamma=1$ is
virtually free.
\end{Theorem}
\begin{proof}
Since $\Gamma$ is finitely presented, we may assume that the
2-skeleton $K^2$ of $K(\Gamma,1)$ is finite.  Let $X$ be the
universal cover of $K^2$ with a lifted metric. By the
\v{S}varc-Milnor Lemma, $X$ is quasi-isometric to $\Gamma$. By
Theorem \ref{asdim1}, $X$ (and hence $\Gamma$) is quasi-isometric to
a tree. By Stallings' theorem $\Gamma$ is virtually free.
\end{proof}
This theorem does not hold for finitely presented groups. For
instance, $\Z_2\wr\Z$ is not virtually free and $\as(Z_2\wr\Z)=1$
\cite{Ge}.

\section{The Hurewicz-type theorem for groups}

The Hurewicz-type theorem (Theorem \ref{as-Hurewicz}) took various
forms in \cite{BD3}:
\begin{Theorem} Let $\Gamma$ be a finitely generated group acting by
isometries on the geodesic metric space $X$. Let $x_0\in X$ and
suppose that for every $R$, the set $\{g\in\Gamma\mid
d(g(x_0),x_0)\le R\}$ has $\as\le n$. Then $\as \Gamma\le \as X+n.$
\end{Theorem}

\begin{Theorem}[Extension Theorem]\label{extension} Let \[\begin{CD}
1 @>>> K @>>> G @>>> H @>>>1\end{CD}\]be an exact sequence with $G$
finitely generated. Then \[\as G\le \as H+\as K\hbox{.}\]
\end{Theorem}

We remark that Dranishnikov and Smith have extended both of these
versions of the Hurewicz-type theorem to all groups. Also they have
shown that for a short exact sequence of abelian groups,
\[\begin{CD} 0 @>>> B @>>> A @>>> C @>>>0\end{CD}\] we have the
equality \[\as A=\as B+\as C\] using Theorem \ref{tensor} below.

The extension theorem above (Theorem \ref{extension}) was used by
Bell and Fujiwara \cite{BellFujiwara} in their upper bound estimates
for mapping class groups of surfaces with genus at most 2:

\begin{Theorem} Let $S_{g,p}$ be an orientable surface of genus $g\le 2$
with $p$ punctures. Suppose that $3g-3+p>1$. Then, \[\as
MCG(S_{g,p})=cd(MCG(S_{g,p}))\hbox{,}\] where $MCG$ denotes the
mapping class group and $cd(\cdot)$ is the cohomological dimension.
\end{Theorem}

Since the braid group $B_n$ is isomorphic to the mapping class group
of a disk with $n$ punctures, we see that a copy of $B_n$ sits
inside $MCG(S_{0,n+1}),$ the mapping class group of the sphere with
$n+1$ punctures. Applying the previous result, we obtain the
following.

\begin{cor} \label{braid} Let $B_n$ be the braid group on $n$ strands.
If $n\ge 3,$ $\as B_n\le n-2$.
\end{cor}

Presently, we discuss other applications of the Extension Theorem.


\section{Polycyclic groups}

\begin{definition} A group $G$ is called {\em polycyclic} \index{polycyclic group}
if there exists a sequence of subgroups $\{1\}=G_0\subset
G_1\subset\dots \subset G_n=G$ such that $G_i\triangleleft G_{i+1}$
and $G_{i+1}/G_i$ is cyclic.

The {\em Hirsch length} \index{Hirsch length} of a polycyclic group,
$h(G)$, is defined to be the number of factors $G_{i+1}/G_i$
isomorphic to $\mathbb{Z}$.
\end{definition}

\begin{theorem} Let $\Gamma$ be a finitely generated polycyclic
group.  Then $\as \Gamma= h(\Gamma).$
\end{theorem}

\begin{proof}  The proof of the inequality $\as \Gamma= h(\Gamma)$ is given
in \cite{BD3}. The inequality $\as\Gamma\ge h(\Gamma)$ was proven in
\cite{DrSm1}. We repeat the first. Denote the sequence of subgroups
satisfying the polycyclic condition by: $\{1\}=\Gamma_0\subset
\Gamma_1\subset\cdots\subset \Gamma_n=\Gamma.$ Then, by Theorem
\ref{extension}, we have

\[\begin{array}{rcl}
\as\Gamma & \le&
\as\Gamma_n/\Gamma_{n-1}+\Gamma_{n-1}\\
&\le& \as\Gamma_n/\Gamma_{n-1}+\as\Gamma_{n-1}/\Gamma_{n-2}+\as
\Gamma_{n-2}\\
&\vdots&\\
&\le &\as\Gamma_n/\Gamma_{n-1}+\cdots+\as\Gamma_{1}/\Gamma_{0}+\as
\Gamma_{0}.
\end{array}\]
Since $\as \Gamma_{i+1}/\Gamma_i$ is only positive if
$\Gamma_{i+1}/\Gamma_i$ is isomorphic to $\mathbb{Z},$ and since in
this case, $\as \Gamma_{i+1}/\Gamma_i=1,$ we conclude $\as \Gamma\le
h(\Gamma).$
\end{proof}

\begin{cor} [\cite{BG}] For polycyclic groups the Hirsch length $h(\Gamma)$ is
a quasi-isometry invariant.
\end{cor}
This result was extended to solvable groups
\cite{Sauer},\cite{Sh}.

Since every finitely generated nilpotent group is polycyclic we
immediately obtain the following result.

\begin{cor} \label{6} Let $\Gamma$ be a finitely generated nilpotent
group.  Then $\as\Gamma= h(\Gamma).$
\end{cor}

More generally, let $G$ be a solvable group with commutator series:
\[1=G_0\subset G_1\subset \cdots\subset G_n=G\hbox{,}\]
so that $G_i=[G_{i+1},G_{i+1}]$. Then, the Hirsch length is defined
to be
\[h(G)=\sum\dim_{\mathbb{Q}}((G_{i+1}/G_i\otimes\mathbb{Q})\hbox{.}\]

\begin{Theorem}\cite[Theorem 3.2]{DrSm1}\label{tensor} For an abelian
group $A$, $\as A=\dim_\mathbb{Q}(A\otimes\mathbb{Q}).$
\end{Theorem}

\begin{cor}\cite[Theorem 3.4]{DrSm1} For a solvable group
$G$, $\as G\le h(G)$.
\end{cor}

\begin{example} Consider the integer Heisenberg group, $H$, i.e. all
matrices of the form \[\begin{pmatrix} 1& x & y\\0 & 1 & z\\0 & 0 &
1\end{pmatrix}\] with the usual multiplication. It is easy to find a
central series for $H$:
\[\begin{pmatrix} 1& 0 & 0\\0 & 1 & 0\\0 & 0 &
1\end{pmatrix}\subset \begin{pmatrix} 1& 0 & 0\\0 & 1 & z\\0 & 0 &
1\end{pmatrix}\subset \begin{pmatrix} 1& x & y\\0 & 1 & z\\0 & 0 &
1\end{pmatrix}.\] So its Hirsch length is $3$ (there are three
copies of $\mathbb{Z}$ in the factor groups). Thus, by Corollary
\ref{6}, $\as H\le 3.$
\end{example}

Modifying the previous example, one could instead take the real
Heisenberg group, $H$. It is the simplest nilpotent Lie group. With
this example in mind, one could extend Corollary \ref{6} to
nilpotent Lie groups $N$ by defining the Hirsch length $h(N)$ as the
sum of the number of factors in $\Gamma_{i+1}/\Gamma_i$ isomorphic
to $\mathbb{R}$ for the central series $\{\Gamma_i\}$ of $N$. We
take an equivariant metric on $N$ and on the quotients. Then the
projection $\Gamma_{i+1}\to \Gamma_{i+1}/\Gamma_i$ is 1-Lipschitz
and $\Gamma_{i+1}/\Gamma_i$ is coarsely isomorphic to
$\mathbb{R}^{n_i}$. Then we have

\begin{cor}  \label{Lie} Let $N$ be a nilpotent Lie group endowed with an equivariant
metric. Then $\as N= h(N).$
\end{cor}

Since $h(N)=\dim N$ for simply connected $N$, we obtain

\begin{cor} \cite[Theorem 3.5]{CG}  \label{cg} For
a simply connected nilpotent Lie group $N$ endowed with an
equivariant metric $\as N= \dim N.$
\end{cor}

Corollary \ref{cg} is the main step in the proof of the following

\begin{theorem}\cite{CG} For a connected Lie group $G$ and its maximal compact subgroup $K$
there is a formula $\as G/K=\dim G/K$ where $G/K$ is endowed with a
G-invariant metric.
\end{theorem}

\section{Groups acting on trees}

\subsection{Bass-Serre Theory}

We briefly discuss the Bass-Serre theory describing the
correspondence between groups acting on trees and generalizations of
amalgamated products. This treatment follows \cite{Se}. For
generalizations of the theory, see \cite{BH} and \cite{Be03}.

Let $Y$ be a nonempty connected graph with vertex set $V(Y)$ and
(directed) edge set $E(Y)$. If $y\in E(Y)$ is an edge, we denote by
$\bar{y}$ the edge $y$ with opposite orientation. The vertex $t(y)$
will denote the terminal vertex of the edge $y$ and the vertex
$i(y)$ will denote the initial vertex of $y.$ We wish to define a
structure called a graph of groups, which can be thought of as a
recipe for building groups in a geometric way. For each $P\in V(Y)$
let $G_P$ be a group and to each $y\in E(Y)$ assign a group $G_{y}$
and two injective homomorphisms $\phi_y:G_y\to G_{t(y)}$ and
$\phi_{\bar{y}}:G_y\to G_{i(y)}.$ Together, the groups,
homomorphisms and graph form the graph of groups \index{graph of
groups} $(G,Y)$.

We want to define a group called the fundamental group of the graph
of groups associated to $(G,Y)$.

First, we define an auxiliary group $F(G,Y)$ associated to the graph
of groups $(G,Y)$. In terms of generators and relations, $F(G,Y)$
can be described as the group generated by all elements of the
vertex groups $G_P$ along with all edges $y\in E(Y)$ of the graph
$Y$. The relations are 1) the relations amongst the groups $G_P$, 2)
the relation that $\bar{y}=y^{-1}$ and 3) an interaction between
edge groups and vertex groups described by
$y\phi_y(a)y^{-1}=\phi_{\bar{y}}(a),$ if $y\in E(Y)$ and $a\in G_y.$

More succinctly, $F(G,Y)$ can be described as the quotient of the
free product $\ast_{P\in V(Y)}G_P\ast \left<y\in E(Y)\right>$ by the
normal subgroup generated by elements of the form
$y\phi_y(a)y^{-1}(\phi_{\bar{y}}(a))^{-1},$ where $y\in E(Y)$ and
$a\in G_y.$

Let $c$ be a path in $Y$ starting at a vertex $P_0$. Let
$y_1,y_2\dots,y_n$ denote the edges of the path in $Y$ with
$t(y_i)=P_i.$ The {\em length} of $c$ is $\ell(c)=n$, its initial
vertex is $i(c)=P_0$ and its terminal vertex is $t(c)=P_n$. A word
of type $c$ in $F(G,Y)$ is a pair $(c,\mu)$ where $c$ is a path as
above and $\mu$ is a sequence $r_0,r_1,\dots,r_n$ with $r_i\in
G_{P_i}$. The associated element of the auxiliary group is
$|c,\mu|=r_0y_1r_1\dots y_nr_n\in F(G,Y)$.

There are two equivalent description of the fundamental group of
$(G,Y)$, but we describe only the one in terms of based loops in
$F(G,Y)$. Let $P_0$ be a fixed vertex of $Y$. The {\em fundamental
group} of $(G,Y)$ is the subgroup $\pi_1(G,Y,P_0)\subset F(G,Y)$
consisting of words associated to loops $c$ in $Y$ based at $P_0$,
i.e., paths with $i(c)=t(c)=P_0$.

\begin{example}
We give three basic examples which we will refer to later.
   \begin{enumerate}
       \item If all vertex groups are trivial, then
       $\pi_1(G,Y,P_0)\cong\pi_1(Y,P_0)$.
       \item Suppose $Y$ is a graph with two vertices $P$ and $Q$
       and one edge $y$ connecting them. Then
       $\pi_1(G,Y,P)=G_P\underset{G_y}{\ast}G_Q.$
       \item Suppose $Y$ is a graph with one vertex and one edge.
       Then $\pi_1(G,Y,P)=G_P\underset{G_y}{\ast}$, the
       HNN-extension.
   \end{enumerate}
\end{example}

Having constructed the fundamental group, $\pi$, of the graph of
groups, we now describe the construction of the Bass-Serre tree
$\tilde{Y}$ on which the group $\pi$ acts by isometries.

Let $T$ be a maximal tree in $Y$ and let $\pi_P$ denote the
canonical image of $G_P$ in $\pi,$ obtained via conjugation by the
path $c,$ where $c$ is the unique path in $T$ from the basepoint
$P_0$ to the vertex $P.$ Similarly, let $\pi_y$ denote the image of
$\phi_y(G_{t(y)})$ in $\pi_{t(y)}.$ Then, set

\[V(\tilde Y)=\coprod_{P\in V(Y)}\pi/\pi_P\] and
\[E(\tilde Y)=\coprod_{y\in E(Y)}\pi/\pi_y.\]

For a more explicit description of the edges, observe that the
vertices $x\pi_{i(y)}$ and $xy\pi_{t(y)}$ are connected by an edge
for all $y\in E(Y)$ and all $x\in\pi.$ Obviously the stabilizer of
the vertices are conjugates of the corresponding vertex groups, and
the stabilizer of the edge connecting $x\pi_{i(y)}$ and
$xy\pi_{t(y)}$ is $xy\pi_yy^{-1}x^{-1},$ a conjugate of the image of
the edge group. This obviously stabilizes the second vertex, and it
stabilizes the first vertex since $y\pi_yy^{-1}= \pi_{\bar y}\subset
\pi_{i(y)}.$ It is known (see \cite{Se}) that the action of left
multiplication on $\tilde Y$ is isometric.

Now we will assume that the graph $Y$ is finite and that the groups
associated to the edges and vertices are finitely generated with
some fixed set of generators chosen for each group.  We let $S$
denote the disjoint union of the generating sets for the groups, and
require that $S=S^{-1}.$  By the norm $\|x\|$ of an element $x\in
G_P$ we mean the minimal number of generators in the fixed
generating set required to present the element $x.$  We endow each
of the groups $G_P$ with the word metric given by $\dist
(x,y)=\|x^{-1}y\|.$  We extend this metric to the group $F(G,Y)$ and
hence to the subgroup $\pi_1(G,Y,P_0)$ in the natural way, by
adjoining to $S$ the collection $\{y, y^{-1}\mid y\in E(Y)\}.$

The principal result in the theory is the following:

\begin{theorem} To every fundamental group of a graph of groups
there corresponds a tree and an action of the fundamental group on
the tree by isometries. To every isometric action of a group on a
tree, there corresponds a graph of groups construction.
\end{theorem}

\subsection{Asdim of amalgams and HNN-extensions}

We want to apply the Hurewicz-type theorem to groups acting on trees
by isometries. To do this, we have to know the structure of the
so-called $R$-stabilizers, i.e., the set
$W_R(x)=\{\gamma\in\Gamma\mid d(\gamma.x,x)\le R\}$. The following
description is straightforward to prove. It appears in \cite{BD2}.

\begin{prop} \label{R-stab Bass-Serre} Let $Y$ be a non-empty, finite, connected graph, and
$(G,Y)$ the associated graph of finitely generated groups.  Let
$P_0$ be a fixed vertex of $Y,$ then under the action of $\pi$ on
$\tilde Y$, the $R$-stabilizer $W_R(1.P_0)$ is precisely the set of
elements of type $c$ in $F(G,Y)$ with $i(c)=P_0,$ and $l(c)\le R.$
\end{prop}

Using the union theorem, we estimate the asymptotic dimension of the
$R$-stabilizers in terms of the asymptotic dimension of the vertex
stabilizers.

\begin{lemma} \label{asdim WR(x)<n} Let $\Gamma$ act on a tree with compact
quotient and finitely generated stabilizers satisfying $\as
\Gamma_x\le n$ for all vertices $x$, then $\as W_R(x_0)\le n$ for
all $R$.
\end{lemma}

Now, by the Hurewicz-type theorem, we have the following theorem.

\begin{Theorem} [\cite{BD2}] Let $\pi$ denote the fundamental group of a finite
graph of groups with finitely generated vertex groups, $V_\sigma$.
Then $$\as\pi\le \max_\sigma\{\as V_\sigma\} +1.$$
\end{Theorem}

This result was extended to complexes of groups in \cite{Be2}.

\begin{example} Let $A$ and $B$ be finitely generated groups with
$\as A\le n$ and $\as B\le n$. Then $\as(A\ast_C B)\le n+1.$
\end{example}

\begin{example} \label{SL(2,Z)} Since ${\rm SL}(2,\mathbb{Z})\cong
\mathbb{Z}_4\ast_{\mathbb{Z}_2}\mathbb{Z}_6$, we see that
$\as(\rm{SL}(2,\mathbb{Z})\le 1.$ Since it is an infinite group,
$\as{\rm SL}(2,\mathbb{Z})=1.$

We knew that $\as {\rm SL}(2,Z)=1$ already, since it is
quasi-isometric to a tree, by the \v{S}varc-Milnor Lemma, but the
point of this example is to apply the estimate for amalgamated
products.
\end{example}

Next, we give an example to show that this upper bound is sharp in
the case of the amalgamated product. We work out the asymptotic
dimension of the free product in a later section.

\begin{example} Using the van Kampen Theorem one can obtain the
fundamental group of the closed orientable surface of genus $2$ as
an amalgamated product of two free groups. Observe that $\as
\pi_1(M_2)=2$, and $\as F=1$ for a free group.
\end{example}

The fact that limit groups have finite asymptotic dimension was
pointed out to the first author by Bestvina to be an easy
consequence of the example above and deep results in the theory of
limit groups. The class of limit groups consists of those groups
which naturally arise in the study of solutions to equations in
finitely generated groups. One definition is the following:  A
finitely presented group $L$ is a {\em limit group} \index{limit
group} if for each finite subset $L_0\subset L$ there is a
homomorphism to a free group which is injective on $L_0.$  For more
information the reader is referred to \cite{BF} and the references
therein.

\begin{prop} Let $L$ be a limit group.  Then $\as L<\infty.$
\end{prop}

\begin{proof} Construct
$L$ (say with height h) via fundamental groups of graphs of groups
where vertices have height ($h-1$) and height 0 groups are free
groups, free abelian groups and surface groups.
%
\end{proof}

Finally, we consider the HNN-extension of a group. Recall that this
corresponds to a fundamental group of a loop of groups.

\begin{example}
Let $A\ast_C$ denote an HNN-extension of the finitely generated
group $A.$ Then $\as A\ast_C\le \as A+1.$
\end{example}

An easy consequence of this example is that one relator groups have
finite asymptotic dimension.

\begin{prop} Let $\Gamma=\left<S\mid
r_1r_2\dots r_n=1\right>$ be a finitely generated group with one
defining relator. Then $\as\Gamma<n+1.$
\end{prop}

\begin{proof} A result of Moldavanskii from 1967 (see \cite{LS} for example)
states that a finitely generated one-relator group is an HNN
extension of a finitely presented group with shorter defining
relator or is cyclic. In the example above we showed that taking
HNN-extensions preserved finite asymptotic dimension and clearly,
all cyclic groups have asymptotic dimension $\le 1$. We conclude
that $\as \Gamma\le n+1.$
\end{proof}
This proof first appeared in \cite{BD4} and then it was rediscovered
in \cite{Mat2}.

\section{A Formula for the asdim of a Free Product} Previously
we saw that $\as A \ast_C B\le \max\{\as A, \as B\}+1$. In the
following example, the inequality is strict (even when $\as C=\as
A$). Let $F_2$ denote a free group on two letters. Then
$\Gamma=F_2\ast_\mathbb{Z}F_2$ is a free group, so $\as\Gamma=1<\as
F_2+1=2$. Here the inclusion map $\mathbb{Z}\to
\left<a,b\mid\right>$ is given by $n\mapsto a^n$.

In this section we focus on free products only; we find an exact
formula for the asymptotic dimension of a free product (amalgamated
over $e$). We begin with two motivating examples.

\begin{example} Since it is finite, $\as \mathbb{Z}_2=0.$ On the
other hand, $\mathbb{Z}_2\ast\mathbb{Z}_2$ is infinite, so
$\as\mathbb{Z}_2\ast \mathbb{Z}_2=1.$
\end{example}

\begin{example} The free group $\mathbb{F}_2$ on two letters has
$\as \mathbb{F}_2=1$, yet it is a free product of two copies of the
asymptotically $1$-dimensional group $\mathbb{Z}$.
\end{example}

%
%
%
%
%

\begin{theorem} Let $A$ and $B$ be finitely generated groups with
$\as A=n$ and $\as B\le n.$  Then, $\as A\ast B =\max\{n,1\}.$
\end{theorem}

Instead of giving the full proof, which can be found in \cite{BDK},
we simply give the ideas behind the proof.

Let $\Gamma$ denote the group $A\ast B$ and let $X$ be the
Bass-Serre tree on which it acts, say by $\pi:\Gamma\to X$ defined
by $\pi(\gamma)=\gamma A$.

Our first observation is that if $W_1$ and $W_2$ are disjoint
bounded sets in $X$ then the sets $\pi^{-1}(W_i)$ are asymptotically
disjoint in $\Gamma$.

Given $\epsilon>0$ we want to construct an $\epsilon$-Lipschitz,
uniformly cobounded map to a uniform $n$-dimensional complex. For a
basepoint $x_0\in X$ we take a cover $\sW=\{W_i\}_i$ of $\Gamma.x_0$
by uniformly bounded disjoint sets. Take these sets so that very
large neighborhoods of these sets have multiplicity $2$.

Each of the sets $\pi^{-1}(W_i)$ has asymptotic dimension $n$, which
is the same as the asymptotic dimension of $A$. Also, any two of
these sets are asymptotically disjoint. Each set $\pi^{-1}(W_i)$
admits an $\epsilon$-Lipschitz, uniformly cobounded map to a
$n$-dimensional complex. Given disjoint sets $W_i$ and $W_j$ whose
large neighborhoods meet, we can find an asymptotic separator for
these sets whose asymptotic dimension is $\le n-1$. Note that here
we use the fact that $\asI X\le \as X$. Thus, there is an
$\epsilon$-Lipschitz uniformly cobounded map from each asymptotic
separator to an $n-1$-complex. We define the mapping from $\Gamma$
to an $n$-dimensional complex by using uniform mapping cylinders
from the inclusion of the asymptotic separator for $W_i$ and $W_j$
into each of $W_i$ and $W_j$. By tweaking the size of the
neighborhoods, one can make the mapping defined this way
$\epsilon$-Lipschitz and it will be uniformly cobounded.

\begin{example} $\as \mathbb{Z}_2\ast\mathbb{Z}_3=1.$
\end{example}
An alternative proof of Theorem 85 is given in \cite{Dr06}.


\section{Coxeter groups}

Let $S$ be a finite set. A Coxeter matrix is a symmetric function
$M:S\times S\to\{1,2,3,\ldots\}\cup\{\infty\}$ with $m(s,s)=1$ and
$m(s,s')=m(s',s)\ge 2$ if $s\neq s'.$ The corresponding Coxeter
group $\sW(M)$ is the group with presentation
\[\sW(M)=\left<S\mid (ss')^{m(s,s')}=1\right>\] where $m(s,s')=\infty$
means no relation. The associated Artin group $\mathcal{A}(M)$ is
the group with presentation
\[\mathcal{A}(M)=\left<S\mid (ss')^{m(s,s')}=(s's)^{m(s,s')}\right>.\]

\begin{Theorem}[\cite{DJ}] Every Coxeter group has finite asymptotic
dimension.\qed
\end{Theorem}
The proof is based on a remarkable embedding theorem of
Januszkiewicz \cite{DJ},\cite{TJan}:

\begin{Theorem} Every Coxeter group $\Gamma$ can be isometrically
embedded in a finite product of trees $\prod T_i$ with the $\ell_1$
metric on it in such a way that the image of $\Gamma$ under this
embedding in contained in the set of vertices of $\prod T_i$.
\end{Theorem}

It should be noted that, except in specific cases, the asymptotic
dimension of a Coxeter group is unknown. It follows from \cite{Dr05}
that the asymptotic dimension of a Coxeter group can be estimated
from below as its virtual cohomological dimension: $\as W(M)\ge
vcd(W(M))$. The upper bound is given by the number of trees in
Januszkiewicz's embedding theorem which equals the number of
generators $|S|$. It was noted in \cite{DSc} that this number can be
lowed to the chromatic number of the nerve $N(M)$ of the Coxeter
group. Recently in \cite{Dr06} it was shown that $\as W(M)\le\dim
N(M)+1$ for even Coxeter groups. A Coxeter group is called {\em
even} if all non-diagonal entries in its Coxeter matrix are even (or
$\infty$).


An Artin group with associated Coxeter matrix $M$,
$\mathcal{A}=\mathcal{A}(M)$, is said to be of finite type if
$\sW(M)$ is finite. It is said to be of affine type if $\sW(M)$ acts
as a proper, cocompact group of isometries on some Euclidean space
with the elements of $S$ acting as affine reflections.

The following approach to finding an upper bound for the $\as$ of
(certain) Artin groups was suggested by Robert Bell.

In \cite{ChCr}, Charney and Crisp observe that each of the Artin
groups $\mathcal{A}(A_n),$ $\mathcal{A}(B_n)$ of finite type and the
Artin groups $\mathcal{A}(\tilde{A}_{n-1})$ and
$\mathcal{A}(\tilde{C}_{n-1})$ of affine type is a central extension
of a finite index subgroup of $\rm{MCG}(S_{0,n+2})$ when $n\ge 3.$
Combining this with the fact that the centers of the Artin groups of
finite type are infinite cyclic and the centers of those of affine
type are trivial gives the following corollary.

\begin{cor}\label{artin-cor} Let $n\ge 3.$ Then if $\mathcal{A}$ is
an Artin group of finite type $A_n$ or $B_n=C_n,$ we have
$\as\mathcal{A}\le n$; if $\mathcal{A}$ is an Artin group of affine
type $\tilde{A}_{n-1}$ or $\tilde{C}_{n-1},$ $\as\mathcal{A}=n-1.$
\end{cor}

This follows easily from the formula for $\as\rm{MCG}(S_{g,p})$ from
\cite{BellFujiwara}.


\section{Hyperbolic groups}

The goal of this section will be to see that a $\delta$-hyperbolic
finitely generated group has finite asymptotic dimension. This
result was announced in Gromov's book \cite{Gr93} and an explicit
proof of more general results appear in \cite{Ro03} and \cite{Ro05}.

\begin{theorem}\label{hyperb} Every finitely generated hyperbolic group has finite
asymptotic dimension.\index{asymptotic dimension!of a hyperbolic
group} \label{hyperbolic group finite asdim}
\end{theorem}

\begin{proof} (cf. Example in Section \ref{asdimT=1})
Let $Y={\rm Cay}(\Gamma, S)$ where $S$ is a symmetric finite
generating set. Then $Y$ is a geodesic metric space with bounded
geometry and $Y$ is quasi-isometric to $\Gamma$ in the word metric
$d_S$. We denote by $e$ the vertex of $Y$ corresponding to the
identity of $\Gamma$.

Let $r\gg \delta$ be given and define concentric annuli
\[A^n=\{x\in Y\mid 2(n-1)r\le d(x,e)\le 2nr\}\] of thickness $2r$
and shells,
\[S^n=\{x\in Y\mid d(x,e)=2nr\}\text{.}\] Let $\{s^n_i\}$ be an enumeration of
the elements of $S^n$.

For $n=1,2$ define $U^n=A^n$. For $n\ge 3$ define families $U^n_i$
as follows:
\[U^n_i=\{x\in A^n\mid \exists \text{ geodesic } [e,x]
\text{ containing } s^{n-2}_i\in S^{n-2}\}.\] It is easy to check
that $\diam(U^n_i)\le 4r$. Next, take $B=B(x,r/2).$ Clearly, $B$ can
meet at most $2$ of the annuli. Suppose $y_i\in U^n_i\cap B$. Then
since $Y$ is $\delta$-hyperbolic, $d(s_i,s)\le \delta$ where $s\in
S^{n-2}$ is on a geodesic $[e,x]$. We conclude that the number of
$U^n_i$ meeting $B$ is bounded by the cardinality number of
geodesics in a $\delta$-ball in $Y$. This number is
$|S|^{\lceil\delta\rceil}$, where $\lceil\delta\rceil$ denotes the
least integer greater than or equal to $\delta$. Thus, the
$r/2$-multiplicity of $\{U^n_i\}$ is no more than
$2|S|^{\lceil\delta\rceil}$. Hence, $\as \Gamma<\infty$.
\end{proof}

We note that this theorem follows from Theorem 47 where the sharp
upper bound for the asymptotic dimension is given. It also follows
from the work of Bonk and Schramm \cite{BS}, who  provide a roughly
quasi-isometric embedding of a hyperbolic metric space with bounded
growth at some scale into a convex subset of hyperbolic space.

A different proof of this theorem was obtained by Buyalo and
Lebedeva \cite{Bu},\cite{BL}. For hyperbolic groups they established
the equality
$$
\as\Gamma=\dim\partial_{\infty}\Gamma+1.
$$
Lebedeva found a formula for asymptotic dimension of product of
hyperbolic groups \cite{Le}
$$
\as(\Gamma_1\times\Gamma_2)=\dim(\partial_{\infty}
\Gamma_1\times\partial_{\infty}\Gamma_2)+2.
$$
Her formula applied to examples of hyperbolic groups with Pontryagin
surfaces as the boundaries \cite{Dr07}  gives an example of Coxeter
hyperbolic groups with
$$\as(\Gamma_1\times\Gamma_2)<\as\Gamma_1+\as\Gamma_2.$$ Metric
spaces satisfying this condition were constructed in \cite{BL}.

Also we remark that Theorem \ref{hyperb} was extended by Bell and
Fujiwara \cite{BellFujiwara} to finite asymptotic dimension of
hyperbolic graphs with a certain type of geodesics.

\section{Relatively hyperbolic groups}

In \cite{Gr87}, Gromov defined relative hyperbolicity. Since then it
has been studied by many authors from many points of view.

The point of this section is an explanation of Osin's work,
\cite{Os}.  The main theorem in that paper is the following.

\begin{Theorem}\label{Osin Thm} Let $G$ be a finitely generated group hyperbolic
with respect to a (finite) collection of subgroups
$\{H_\lambda\}_{\lambda\in\Lambda}$. Suppose that each of the groups
$H_\lambda$ has finite asymptotic dimension. Then $\as G<\infty$.
\end{Theorem}

There are several ways of defining relatively hyperbolic groups. We
give the definition that Osin uses in his paper, see also
\cite{bowditch:relatively,Farb,Gr87}.

Let $G$ be a group, $\{H_\lambda\}_{\lambda\in\Lambda}$ a collection
of subgroups of $G$ and $X$ a subset of $G$. We say that $X$ is a
{\em relative generating set of $G$ with respect to}
$\{H_\lambda\}_{\lambda\in\Lambda}$ if $G$ is generated by $X$ and
the union of the $H_\lambda$.

Let $F=(\ast_{\lambda\in\Lambda}H_\lambda)\ast F(X)$, where $F(X)$
is the free group on the set $X$. A {\em relative presentation} for
$G$ is a presentation of the form
\[\left<X,H_\lambda,\lambda\in\Lambda\mid \mathcal{R}\right>\hbox{.}\]
We say that this presentation is {\em finite} and say that $G$ is
{\em finitely presented relative to the collection of subgroups
$\{H_\lambda\}_{\lambda\in\Lambda}$} if $\sharp X$ and $\sharp
\mathcal{R}$ are finite.

Let $\sH=\bigsqcup_{\lambda\in\Lambda} (H_\lambda\setminus\{1\})$.
Given a word $W$ in the alphabet $X\cup\sH$ such that $W$ represents
$1$ in $G$, there is an expression
\[W=_F\prod_{i=1}^kf^{-1}_iR^{\pm1}_if_i\] (with equality in the
group $F$), where $R_i\in\mathcal{R}$ is a relation, and $f_i\in F$
for $i=1,\ldots,k.$ The smallest possible $k$ in such a presentation
is called the relative area of $W$ and is denoted by
$Area^{rel}(W)$.
\begin{defn} A group $G$ is said to be {\em hyperbolic relative to a
collection of subgroups $\{H_\lambda\}_{\lambda\in\Lambda}$} if $G$
is finitely presented relative to
$\{H_\lambda\}_{\lambda\in\Lambda}$  and there is a constant $L>0$
such that for any word $W$ in $X\cup\sH$ representing the identity
in $G$, we have $Area^{rel}(W)\le L\|W\|$.
\end{defn}

Osin showed in \cite{Os1} that when $G$ is generated by a finite set
in the ordinary sense and finitely presented relative to a
collection of subgroups $\{H_\lambda\}_{\lambda\in\Lambda}$ then
$\Lambda$ is known to be finite and the subgroups $H_\lambda$ are
known to be finitely generated.

To prove Theorem \ref{Osin Thm}, Osin applies the Hurewicz-type
theorem. Let $\Gamma(G,X\cup\sH)$ denote the Cayley graph of $G$
with respect to the generating set $X\cup\sH$. Note that
$\Gamma(G,X\cup\sH)$ is hyperbolic, see \cite[Theorem 1.7]{Os1}. It
is not locally finite, but on the other hand, Osin proves that it
has finite asymptotic dimension. Next, he proves that the asymptotic
dimension of so-called relative balls $B(n)$ does not exceed the
maximum of the asymptotic dimensions of the $H_\lambda$. Here, by a
relative ball we mean a set of the form
\[B(n)=\{g\in G\mid |g|_{X\cup\sH}\le n\}\hbox{.}\] In other words
it is the ball of radius $n$ in $G$ with respect to the distance
$d_{X\cup\sH}$ centered at $1$.

\begin{proof}[Proof of Theorem \ref{Osin Thm}] The group $G$ acts
on $\Gamma(G,X\cup\sH)$ by left multiplication. The $R$-stabilizer
at $1$ coincides with the relative ball $B(R)$. By the Hurewicz-type
theorem for groups, $\as G<\infty$.
\end{proof}

A notion related to relative hyperbolicity is that of weak relative
hyperbolicity. A group $G$ is said to be {\em weakly relatively
hyperbolic} with respect to a collection
$\{H_\lambda\}_{\lambda\in\Lambda}$ of subgroups if the Cayley graph
$\Gamma(G,X\cup\sH)$ is hyperbolic, where $X$ is a finite generating
set for $G$ modulo $\{H\}_{\lambda\in\Lambda}$ and $\sH$ is defined
as above.

Although relative hyperbolicity implies weak relative hyperbolicity,
the converse does not hold. The natural question that arises is
whether weak relative hyperbolicity is enough for finite asymptotic
dimension in the sense of Theorem \ref{Osin Thm}. Osin answers this
question in the negative. In particular he proves the following
theorem.

\begin{Theorem} There exists a finitely presented boundedly
generated group of infinite asymptotic dimension.
\end{Theorem}

Recall that a group is said to be {\em boundedly generated} if there
are elements $x_1,\ldots,x_n$ of $G$ such that any $g\in G$ can be
represented in the form $g=x_1^{\alpha_1}\cdots x_n^{\alpha_n}$ for
some $\alpha_1,\ldots,\alpha_n\in\mathbb{Z}$. So with respect to any
generating set $X$ and $\{H_\lambda\}=\left<x_\lambda\right>$ we see
that the Cayley graph $\Gamma(G,X\cup\sH)$ has finite diameter. Thus
as a corollary, we obtain the following.

\begin{Corollary} There exists a finitely presented group $G$ and a
finite collection of subgroups $\{H_\lambda\}_{\lambda\in\Lambda}$
such that
   \begin{enumerate}
       \item each $H_\lambda$ is cyclic (so $\as H_\lambda\le 1$);
       \item the Cayley graph $\Gamma(G,X\cup\sH)$ has finite
diameter -- so it is hyperbolic and $G$ is weakly relatively
hyperbolic with respect to $\{H_\lambda\}$; and
       \item $\as G=\infty$.
   \end{enumerate}
\end{Corollary}

\section{Arithmetic groups}

Let $\Gamma$ be an arithmetic subgroup of a linear algebraic group
$\mathbf{G}$ defined over $\mathbb{Q}$. L. Ji \cite{Ji} generalized
the result of Carlsson and Goldfarb, (Theorem \ref{cg}) to show that
$\Gamma$ has finite asymptotic dimension.

We begin by recalling Carlsson and Goldfarb's theorem.

\begin{theorem} Let $G$ be a connected Lie group with maximal
compact subgroup $K$. Let $X=G/K$ be the associated homogeneous
space endowed with a $G$-invariant Riemannian metric. Then $\as
X=\dim X$.
\end{theorem}

As a corollary, Ji observes the following result.

\begin{Corollary} Let $G$ be a connected Lie group and
$\Gamma\subset G$ any finitely generated discrete subgroup. Then
$\as \Gamma<\infty$.
\end{Corollary}

\begin{proof} The group $\Gamma$ acts properly isometrically on the
proper metric space $X=G/K$. Observe that the map $\gamma\mapsto
\gamma.x_0$ is a coarse equivalence between $\Gamma$ and
$\Gamma.x_0\subset X$ for any choice of $x_0\in X$, (cf. the
\v{S}varc-Milnor Lemma, \cite{dlH}). Thus $\as\Gamma\le \as X$.
\end{proof}

Let $\mathbf{G}$ be a linear algebraic group defined over
$\mathbb{Q}$, i.e., an algebraic subgroup of $GL(n,\mathbb{C})$
defined by polynomial equations with rational coefficients. A
subgroup $\Gamma$ of $\mathbf{G}(\mathbb{Q})$ is called {\em
arithmetic} if $\Gamma$ is commensurable with $\mathbf{G}\cap
GL(n,\mathbb{Z})$.

The following finite dimensionality result essentially follows from
the previous corollary. In addition, Ji gives a lower bound for the
asymptotic dimension that can be achieved when the lattice is
non-cocompact. Let $\rho$ denote the $\mathbb{Q}$-rank of
$\mathbf{G}$, i.e., the maximal dimension of $\mathbb{Q}$-split tori
in $\mathbf{G}$.

\begin{theorem} Let $\mathbf{G}$ be a connected linear algebraic
group defined over $\mathbb{Q}$. Let $\Gamma$ be an arithmetic
subgroup of $\mathbf{G}$, which is assumed to be a lattice in
$G=\mathbf{G}(\mathbb{R})$. Then
\[\dim X-\rho\le \as \Gamma\le \dim X\hbox{.}\]
\end{theorem}

Kleiner observed that the lower bound follows from a result of Borel
and Serre \cite{BorelSerre} and the fact that $\as \Gamma\ge
cd(\Gamma)$, where $cd(\cdot)$ denotes the cohomological dimension,
see \cite{Gr93}.

Later Lizhen Ji extended his result to $S$-arithmetic groups
\cite{Ji2}:
\begin{theorem}
$\as\Gamma<\infty$ for all $S$-arithmetic groups $\Gamma$.
\end{theorem}

\section{Buildings}

In his thesis, D. Matsnev \cite{Mat1} applied the Hurewicz-type
theorem for asymptotic dimension to show that affine buildings
have finite asymptotic dimension.  The method is to reduce the
problem to computing the asymptotic dimension of a certain matrix
group that is coarsely equivalent to a given affine building.

Throughout this section we let $K$ field with a discrete valuation
$\nu$, i.e., a surjection $\nu:K^*\to \mathbb{Z}$ satisfying
\[\nu(x+y)\ge \min\{\nu(x),\nu(y)\}\] for all $x,y\in K^*$.
Let $X$ denote the building associated to $SL(n,K)$, see
\cite[Chapter V.8]{Brown} for a description of $X$.

We now describe the notion of distance on $SL(n,K)$. The ``metric''
on $SL(n,K)$ is actually not a metric at all, but rather a
pseudometric. It inherits this pseudometric from a length function
$\ell$ defined in terms of the discrete valuation $\nu$ as follows:
\[\ell(g)=-\min_{1\le i,j\le n}\{\nu(g_{ij}),\nu(g^{ij})\}\hbox{.}\]
Here $g_{ij}$ is the $ij$-th entry of $g$ and $g^{ij}$ is the
$ij$-th entry of its inverse. The psuedometric is then defined to be
\[\dist(g,h)=\ell(g^{-1}h)\hbox{.}\]

The following is a sketch of the proof of Matsnev's theorem:

\begin{theorem} The affine building $X$ has finite asymptotic
dimension.
\end{theorem}

The affine building $X$ is coarsely equivalent to the group
$G=SL(n,K)$, with the pseudometric described above so it suffices to
compute $\as G$. If $C$ denotes the maximal compact subgroup of $G$,
then we can write $G=CB$, where $B$ is the subgroup of upper
triangular matrices. Since $C$ is compact, $G$ is coarsely
equivalent to $B$, so it remains to compute $\as B$.

To this end, Matsnev defines a map $f:B\to A$ where $A$ denotes the
diagonal matrices. This map simply takes an upper-triangular matrix
to the diagonal matrix obtained by replacing all off-diagonal
entries by $0$. It can be shown that this map is $1$-Lipschitz and
that the set $f^{-1}(B_R(a))$ has asymptotic dimension $0$ uniformly
(in $a$) for every $R$. Thus, by the Hurewicz-type theorem, we see
that $\as B\le \as A$. It is shown that $A$ in the pseudometric is
coarsely equivalent to a subgroup isomorphic to $\mathbb{Z}^{n-1}$.
Although the metric on $\mathbb{Z}^{n-1}$ is not the standard one,
it is Lipschitz equivalent to the standard one and so we conclude
that $\as A=n-1$. Since $A\subset B$, we conclude that $\as B=n-1$
and therefore that $\as X=n-1$.

Recently Jan Dymara and Thomas Schick extended Matsnev's result to
general (Tits) buildings \cite{DySc}.
\begin{theorem}
The asymptotic dimension of a building $X$ equals the asymptotic
dimension of the apartment.
\end{theorem}

\section{Infinite dimensional groups}
It is not at all difficult to find examples of finitely generated
groups with infinite asymptotic dimension. Indeed, let $\Gamma$ be
a finitely generated group containing an isomorphic copy of
$\mathbb{Z}^m$ for each $m$, then $\as Y=\infty$. Some attempts
were made to start up a theory of asymptotically infinite
dimensional spaces by analogy with topology. In \cite{Dr1} an
asymptotic property C was defined as follows: A metric space $X$
has asymptotic property C if for any sequence of natural numbers
$n_1<n_2<\dots$ there is a finite sequence of uniformly bounded
families $\{\sU_i\}^n_{i=1}$ such that the union
$\cup^n_{i=1}\sU_i$ is a cover of $X$ and each $\sU_i$ is
$n_i$-disjoint. T. Radul introduced a notion of transfinite
asymptotic dimension $trasdim$ and proved that $X$ has property C
if and only if $trasdim(X)$ is defined \cite{Ra}. In view of
Borst's theorem which states that the transfinite dimension is
defined for weakly infinite dimensional compacta and his recent
example \cite{Borst} this result shows a striking difference
between asymptotic and topological dimension for infinite
dimensional spaces.

Another way to deal with asymptotically infinite dimensional
groups is to study the dimension function $ad_X(\lambda)$ of a
metric space $X$. We define this function as follows:

\begin{defn} Let $X$ be a metric space and $\sU$ a cover of $X$.
Denote the multiplicity of $\sU$ by $m(\sU)$, i.e.,
$m(\sU)=\sup\card\{U\in\sU\mid x\in U\}$. Denote the Lebesgue number
of $\sU$ by $L(\sU)$, i.e., $L(\sU)$ is the largest number so that
for any $A\subset X$ with $\diam(A)\le L(\sU)$ there is a $U\in\sU$
with $A\subset U$. Define the dimension function of $X$ by:
\[ad_X(\lambda)=\min\{m(\sU)\mid L(\sU)\ge\lambda\}-1\hbox{.}\]
\end{defn}

It is easy to see that $ad_X(\lambda)$ is monotone and
$\lim_{\lambda\to\infty} ad_X(\lambda)=\as X$.

In the section on coarse embeddings we proved that a metric space
with finite asymptotic dimension admits a coarse embedding into
Hilbert space. Although we proved this result directly, in that
discussion we mentioned that what is really needed to prove
embeddability into Hilbert space is Yu's Property A. This property
can be defined in terms of anti-{\v C}ech approximations. Namely a
metric space has Property A if it admits an anti-{\v C}ech
approximation $\{\sU_i\}$ such that the canonical projections to the
nerves $p_{\sU_i}$ are $\epsilon_i$-Lipschitz with $\epsilon_i\to 0$
where the nerves $Nerve(\sU_i)$ are taken with the induced metric
from $\ell_1(\sU_i)$ \cite{Dr-book}.

In \cite{Dr1},\cite{Dr04}, Dranishnikov proved the following
generalization of this embedding result.

\begin{Theorem} A metric space $X$ has property A in each of the
following cases:

(a) $X$ is a discrete metric space with polynomial dimension
growth;

(b) $X$ has property $C$.
\end{Theorem}

\begin{Corollary} Let $\Gamma$ be a finitely generated group whose
dimension function grows polynomially. Then the Novikov conjecture
holds for $\Gamma$.
\end{Corollary}

In order for these results to be interesting, we need examples of
spaces with polynomial dimension growth. Our example (following Roe
\cite{Ro03}) is the restricted wreath product of $\mathbb{Z}$ by
$\mathbb{Z}$.

Let $H$ be the set of finitely supported maps
$\mathbb{Z}\to\mathbb{Z}$. Let $u$ and $v$ be the permutations of
$H$ defined by \[uf(n)=f(n)+\delta_{n0},\qquad
vf(n)=f(n+1)\hbox{.}\] Let $G$ be the group generated by $u$ and
$v$. Endow $G$ with the word metric, i.e., define $d(g,h)$ to be the
length of the shortest presentation of the element $g^{-1}h$ in the
alphabet $\set{u,v}$. It is easy to see that the elements
\[u,v^{-1}uv,\ldots,v^{(-n-1)}u,v^{n-1}\]
generate an isomorphic copy of $\mathbb{Z}^n$ for each $n$, so $\as
G=\infty$. On the other hand, Dranishnikov proved the following
theorem in \cite{Dr04}.

\begin{Theorem} Let $N$ be a finitely generated nilpotent group and
let $G$ be a finitely generated group with $\as G<\infty$. Then the
restricted wreath product $N\wr G$ has polynomial dimension growth.
\end{Theorem}

Another famous group with infinite asymptotic dimension is
Thompson's group $F$. It has many incarnations, but the easiest
description (combinatorially) is that it is the group whose
presentation is \[F=\left<x_0,x_1,x_2,\ldots\mid
x_{j+1}=x_i^{-1}x_jx_i\hbox{, for $i<j$}\right>\hbox{.}\] Notice
that for $i\ge 2$, $x_i=x_0^{1-i}x_1x^{i-1}$, so $F$ is finitely
generated. The growth rate of the dimension function of Thompson's
group F is unknown. The infinite dimensionality of the Thompson
group is based on the fact that $F$ contains $\Z^n$ as a subgroup
for all $n$. An example of a torsion group with infinite asymptotic
dimension is Grigorchuk's group \cite{Sm2}. We briefly recall a fact
from \cite{Dr04} about the growth of the dimension function as it
applies to groups.

We saw that a particular wreath product had polynomial dimension
growth. Gromov's group \cite{Gr03} containing an expander has
exponential dimension growth. The next proposition says that this is
the fastest the function can grow.

\begin{prop} \cite[Proposition 2.1]{Dr04} Let $\Gamma$ be a finitely generated group. Then there
is an $a>0$ so that $ad_{\Gamma}(\lambda)\le e^{a\lambda}.$
\end{prop}

\begin{proof} There is an $a>0$ for which $|B_\lambda(x)|\le
e^{a\lambda}$. Clearly the cover of $\Gamma$ by all balls
$B_{\lambda}(x)$ has Lebesgue number at least $\lambda$. Also, the
sets are uniformly bounded and have multiplicity $\le
|B_{\lambda}(x)|\le e^{a\lambda}$, as required.
\end{proof}

Also, the growth of the asymptotic dimension function is not a
coarse invariant, but it is an invariant of quasi-isometries. Since
any two word metrics on a group are quasi-isometric, the growth of
the dimension function is a group invariant.

\def\polhk#1{\setbox0=\hbox{#1}{\ooalign{\hidewidth
 \lower1.5ex\hbox{`}\hidewidth\crcr\unhbox0}}}
\providecommand{\bysame}{\leavevmode\hbox
to3em{\hrulefill}\thinspace}
\providecommand{\MR}{\relax\ifhmode\unskip\space\fi MR }
\providecommand{\MRhref}[2]{%
 \href{http://www.ams.org/mathscinet-getitem?mr=#1}{#2}
} \providecommand{\href}[2]{#2}

\end{document}